\documentclass[a4paper,12pt]{amsart}
\usepackage{amssymb, amsfonts, amsthm}
\usepackage{ifthen}
\usepackage{graphicx}
\nonstopmode \numberwithin{equation}{section}
\usepackage[usenames]{color}
\newtheorem{thm}{Theorem}%[section]
\newtheorem{cor}{Corollary}%[section]%[equation]
\newtheorem{lem}{Lemma}%[section]

\newtheorem{rem}{Remark}%[section][equation]
\newtheorem{rems}[equation]{Remarks}
\theoremstyle{definition}
\newtheorem{defin}{Definition}%[section]%[equation]
\newtheorem{examp}[equation]{Example}%[section]
\newtheorem{prob}[equation]{Problem}
\newtheorem{ques}[equation]{Question}
\newtheorem{op}{Open Problem}%[section]%[equation]

\newtheorem{conj}[equation]{Conjecture}
\newtheorem{deter}[equation]{Determination}
\newtheorem{case}{Case}[section]%[equation]
\newtheorem{subcase}[equation]{Subcase}
\newtheorem{claim}{Claim}[section]%[equation]
\newtheorem{subclaim}{Subclaim}

\newcounter {own}
\def\theown {\thesection       .\arabic{own}}

\newenvironment{pf}[1][]{%
 \vskip 3mm
 \noindent
 \ifthenelse{\equal{#1}{}}%
  {{\bf Proof. }}%
  {{\bf #1.} }%
 }%
{\qed\bigskip}

\newcounter{alphabet}
\newcounter{tmp}
\newenvironment{Thm}[1][]{\refstepcounter{alphabet}%
\bigskip%
\noindent%
{\bf Theorem \Alph{alphabet}}%
\ifthenelse{\equal{#1}{}}{}{ (#1)}%
{\bf .} \itshape}{\vskip 8pt}
\makeatletter
\newcommand{\Ref}[1]{\@ifundefined{r@#1}{}{\setcounter{tmp}{\ref{#1}}\Alph{tmp}}}
\makeatother

\newenvironment{Lem}[1][]{\refstepcounter{alphabet}%
\bigskip%
\noindent%
{\bf Lemma \Alph{alphabet}}%
{\bf .} \itshape}{\vskip 8pt}

\newcommand{\IR}{{\mathbb R}}

\newcommand{\diam}{{\operatorname{diam}}}

\def\be{\begin{equation}}
\def\ee{\end{equation}}

\newcommand{\bee}{\begin{enumerate}}
\newcommand{\eee}{\end{enumerate}}

\newcommand{\blem}{\begin{lem}}
\newcommand{\elem}{\end{lem}}
\newcommand{\bthm}{\begin{thm}}
\newcommand{\ethm}{\end{thm}}
\newcommand{\bcor}{\begin{cor}}
\newcommand{\ecor}{\end{cor}}
\newcommand{\beg}{\begin{examp}}
\newcommand{\eeg}{\end{examp}}
\newcommand{\begs}{\begin{examples}}
\newcommand{\eegs}{\end{examples}}
\newcommand{\bdefe}{\begin{defin}}
\newcommand{\edefe}{\end{defin}}
\newcommand{\bprob}{\begin{prob}}
\newcommand{\eprob}{\end{prob}}
\newcommand{\bques}{\begin{ques}}
\newcommand{\eques}{\end{ques}}
\newcommand{\bei}{\begin{itemize}}
\newcommand{\eei}{\end{itemize}}

\newcommand{\bde}{\begin{deter}}
\newcommand{\ede}{\end{deter}}
\newcommand{\bca}{\begin{case}}
\newcommand{\eca}{\end{case}}
\newcommand{\bsca}{\begin{subcase}}
\newcommand{\esca}{\end{subcase}}
\newcommand{\bcl}{\begin{claim}}
\newcommand{\ecl}{\end{claim}}
\newcommand{\bscl}{\begin{subclaim}}
\newcommand{\escl}{\end{subclaim}}
\newcommand{\bcon}{\begin{conj}}
\newcommand{\econ}{\end{conj}}
\newcommand{\bcons}{\begin{conjs}}
\newcommand{\econs}{\end{conjs}}
\newcommand{\bprop}{\begin{propo}}
\newcommand{\eprop}{\end{propo}}
\newcommand{\br}{\begin{rem}}
\newcommand{\er}{\end{rem}}
\newcommand{\brs}{\begin{rems}}
\newcommand{\ers}{\end{rems}}
\newcommand{\bo}{\begin{obser}}
\newcommand{\eo}{\end{obser}}
\newcommand{\bos}{\begin{obsers}}
\newcommand{\eos}{\end{obsers}}
\newcommand{\bpf}{\begin{pf}}
\newcommand{\epf}{\end{pf}}
\newcommand{\ba}{\begin{array}}
\newcommand{\ea}{\end{array}}
\newcommand{\beq}{\begin{eqnarray}}
\newcommand{\beqq}{\begin{eqnarray*}}
\newcommand{\eeq}{\end{eqnarray}}
\newcommand{\eeqq}{\end{eqnarray*}}

\newcommand{\bop}{\begin{op}}
\newcommand{\eop}{\end{op}}

\newtheorem{pfofThm1.5}[equation]{}

\newcounter{minutes}
\divide\time by 60
\newcounter{hours}
\multiply\time by 60 \addtocounter{minutes}{-\time}
\begin{document}

\bibliographystyle{amsplain}

\title{On the subinvariance of uniform domains in metric spaces}

%%%%%%%% BEGIN TIMESTAMP
\def\thefootnote{}
\footnotetext{ \texttt{\tiny File:~\jobname .tex,
          printed: \number\year-\number\month-\number\day,
          \thehours.\ifnum\theminutes<10{0}\fi\theminutes}
} \makeatletter\def\thefootnote{\@arabic\c@footnote}\makeatother
%%%%%%%% END TIMESTAMP

\author{Yaxiang  Li}
\address{Yaxiang Li,  College of Science,
Central South University of
Forestry and Technology, Changsha,  Hunan 410004, People's Republic
of China} \email{yaxiangli@163.com}

\author{Manzi Huang}
\address{Manzi Huang, Department of Mathematics,
Hunan Normal University, Changsha,  Hunan 410081, People's Republic
of China} \email{mzhuang79@yahoo.com.cn}

\author{Xiantao Wang${}^{\mathbf{*}}$}
\address{Xiantao Wang, Department of Mathematics,
Shantou University, Shantou, Guangdong 515063, People's Republic
of China} \email{xtwang@stu.edu.cn}

\author{Qingshan Zhou}
\address{Qingshan Zhou, Department of Mathematics,
Hunan Normal University, Changsha,  Hunan 410081, People's Republic
of China} \email{q476308142@qq.com}

\date{}
\keywords{Subinvariance, uniform domain, weak quasisymmetry, free quasiconformality, quasihyperbolic homemorphism, locally John domain.
\\
${}^{\mathbf{*}}$ Corresponding author}

\begin{abstract}
Suppose that $X$ and $Y$ are quasiconvex and complete metric spaces, that $G\subset X$ and $G'\subset Y$ are
domains, and that $f: G\to G'$ is a homeomorphism. Our main result is the following
subinvariance property of the class of uniform domains: Suppose both $f$ and $f^{-1}$ are weakly quasisymmetric mappings and $G'$ is
a quasiconvex domain. Then the image $f(D)$ of every uniform subdomain $D$ in $G$ under $f$ is uniform.
The subinvariance of uniform domains with respect to freely quasiconformal mappings or quasihyperbolic mappings is also studied with the additional condition that both $G$ and $G'$ are locally John domains.
\end{abstract}

\thanks{The research was partly supported by  NNSF of China}

\maketitle\pagestyle{myheadings} \markboth{Yaxiang Li, Manzi Huang, Xiantao Wang and Qingshan Zhou}{On the subinvariance of uniform domains in metric spaces}

\section{Introduction and main results}\label{sec-1}

The quasihyperbolic metric (briefly, QH metric) was introduced by Gehring and
his students Palka and Osgood in the 1970's \cite{Geo, GP} in the setting of Euclidean spaces
${\mathbb R}^n$ $(n\ge 2).$ Since its first appearance, the quasihyperbolic metric has become
an important tool in the geometric function theory of Euclidean spaces, especially, in the study of quasiconformal and quasisymmetric mappings. From late 1980's onwards, V\"ais\"al\"a has developed the theory of (dimension) free quasiconformal mappings (briefly, free theory) in Banach spaces, which is based on properties of the quasihyperbolic metric \cite{Vai6-0, Vai6, Vai6', Vai7, Vai8}. The main advantage of this approach over generalizations based on the conformal modulus (see \cite{Heinonen-book} and references therein) is that it does not make use of volume integrals, thus allowing the study of quasiconformality in infinite dimensional Banach spaces.

The class of quasisymmetric mappings on the real axis was first introduced by Beurling and
Ahlfors \cite{BA}, who found a way to obtain a quasiconformal extension of a quasisymmetric self-mapping of the
real axis to a self-mapping of the upper half-plane. This idea
was later generalized by Tukia and V\"ais\"al\"a, who studied quasisymmetric
mappings between metric spaces \cite{TV}. In
1998, Heinonen and Koskela \cite{HK} proved a remarkable result, showing that the concepts of quasiconformality and
quasisymmetry are quantitatively equivalent in a large class of metric spaces, which
includes Euclidean spaces.  Also, V\"ais\"al\"a proved the quantitative equivalence between free quasiconformality and
quasisymmetry of homeomorphisms between two Banach spaces. See \cite[Theorem 7.15]{Vai8}. Against this background, it is not surprising that the study of quasisymmetry in metric spaces has recently attracted significant attention.

The main goal of this paper is to establish the subinvariance of uniform domains in suitable metric spaces with respect to quaisymmetric mappings, freely quasiconformal mappings and quasihyperbloic mappings. We start by recalling some basic definitions.
Through this paper, we always assume that $X$ and $Y$ are metric spaces. We follow the notations and terminology of \cite{HK-1, HK, HL, Tys, Vai8}.

\bdefe\label{japan-30}
A homeomorphism $f$ from $X$ to $Y$ is said to be \begin{enumerate}
\item
{\it quasiconformal} if there is a constant $H<\infty$ such that
\be\label{mon-1}\limsup_{r\to 0}\frac{L_f(x,r)}{l_f(x,r)}\leq H\ee
for all $x\in X$;
\item
{\it quasisymmetric} if there is a constant $H<\infty$ such that
\be\label{mon-2}\frac{L_f(x,r)}{l_f(x,r)}\leq H\ee
for all $x\in X$ and all $r>0$,\end{enumerate}
 where $L_f(x,r)=\sup_{|y-x|\leq r}\{|f(y)-f(x)|\}$
and
$l_f(x,r)=\inf_{|y-x|\geq r}\{|f(y)-f(x)|\}$.
\edefe
\noindent Here and in what follows, we always use $|x-y|$ to denote the distance between $x$ and $y$.

\bdefe\label{japan-31} A homeomorphism $f$ from $X$ to $Y$ is said to be
\begin{enumerate}
\item $\eta$-{\it quasisymmetric} if there is a homeomorphism $\eta : [0,\infty) \to [0,\infty)$ such that
$$ |x-a|\leq t|x-b|\;\; \mbox{implies}\;\;   |f(x)-f(a)| \leq \eta(t)|f(x)-f(b)|$$
for each $t\geq 0$ and for each triple $x,$ $a$, $b$ of points in $X$;

\item {\it weakly $H$-quasisymmetric} if
$$ |x-a|\leq |x-b|\;\;  \mbox{ implies}\;\;   |f(x)-f(a)| \leq H|f(x)-f(b)|$$
for each triple $x$, $a$, $b$ of points in $X$.
\end{enumerate}\edefe
%%%%%%%%%%%%%%%%%

\br\label{japan-32} The following observations follow immediately from Definitions \ref{japan-30} and \ref{japan-31}.
\begin{enumerate}
\item\label{hwz-1}
The quasisymmetry implies the quasiconformality;
\item\label{hwz-2}
A homeomorphism $f$ from $X$ to $Y$ is quasisymmetric with coefficient $H$ defined by Definition \ref{japan-30}\eqref{mon-2} if and only if it is weakly $H$-quasisymmetric;
\item\label{hwz-3}
The $\eta$-quasisymmetry implies the weak $H$-quasisymmetry with $H=\eta(1)$. Obviously, $\eta(1)\geq 1$. In general, the converse is not true (cf. \cite[Theorem $8.5$]{Vai8}). See also \cite{hprw} for the related discussions.
\end{enumerate}\er

%%%%%%%%%%%%%%%%%
%%%%%%%%%%%%%%%%%

The definition of free quasiconformality is as follows.

\bdefe\label{japan-33} Let $G\varsubsetneq X$ and $G'\varsubsetneq Y$ be two domains (open and connected), and let $\varphi:[0,\infty)\to [0,\infty)$ be a homeomorphism with $\varphi(t)\geq t$. We say
that a homeomorphism $f: G\to G'$ is \begin{enumerate}
\item \label{sunday-1}
{\it $\varphi$-semisolid } if $$k_{G'}(f(x),f(y))\leq \varphi(k_G(x,y))$$
for all $x$, $y\in G$;

\item \label{sunday-2} $\varphi$-{\it solid} if both $f$ and ${\it f^{-1}}$ are $\varphi$-semisolid;

\item {\it freely
$\varphi$-quasiconformal} ($\varphi$-FQC in brief) or {\it fully $\varphi$-solid}
if $f$ is
$\varphi$-solid in every subdomain of $G$,\end{enumerate}
where $k_G(x,y)$ denotes the quasihyperbolic distance of $x$ and $y$ in $G$. See Section \ref{sec-2} for the precise definitions of $k_G(x,y)$ and other
notations and concepts in the rest of this section.
\edefe

It follows from \cite[Remark, p. 121]{FHM} and \cite[Theorem 5.6]{Vai2} that uniform domains are subinvariant with respect to quasiconformal mappings in $\IR^n$ ($n\geq 2$). By this, we mean that if $f:\; G\to G'$ is a $K$-quasiconformal mapping, where $G$ and $G'$ are domains in $\IR^n$, and if $G'$ is $c$-uniform,
then $D'=f(D)$ is $c'$-uniform
 for every $c$-uniform subdomain
 $D \subset G$,
 where $c'=c'(c, K, n)$ which means that the constant $c'$ depends only on the coefficient  $c$ of the uniformity of $D$, the coefficient $K$ of quasiconformality of $f$ and the dimension $n$ of the Euclidean space $\IR^n$.
In the free theory, V\"ais\"al\"a also studied this property of uniform domains in Banach spaces and proved the following result.

\begin{Thm}\label{ThmB}  $($\cite[Theorem 2.44]{Vai4}$)$\quad
Suppose that $G\subset E$ and $G'\subset E'$, where $E$ and $E'$ are Banach spaces with dimension at least $2$,
that the domains $G'$ and $D\subset G$ are $c$-uniform, and that $f:\; G\to G'$ is $M$-quasihyperbolic.
Then $f(D)$ is $c'$-uniform with $c'=c'(c,M)$.
\end{Thm}

In 2012, Huang, Vuorinen and Wang proved the
subinvariance property of uniform domains is also true
with respect to freely quasiconformal mappings (\cite[Theorem 1]{HVW}).
See \cite{BHX, GM, Vai2, Vai8, Xie} for similar discussions in this line.

Our work is motivated by the above ideas which we will extend to the
context of weakly quasisymmetric mappings, freely quasiconformal mappings and quasihyperbolic mappings in metric spaces. Our first result is as follows.

\begin{thm}\label{thm1.1}
Suppose that $X$ and $Y$ are quasiconvex and complete metric spaces, that $G\varsubsetneq X$ is a domain,
$G'\varsubsetneq Y$ is a quasiconvex domain, and that both $f: G\to G'$ and $f^{-1}: G'\to G$ are weakly quasisymmetric mappings.
For each subdomain $D$ of $G$, if $D$ is uniform, then $D'=f(D)$ is uniform, where the coefficient of uniformity of $D'$ depends
only on the given data of $X$, $Y$, $G$, $G'$, $D$, $f$ and $f^{-1}$.
\end{thm}

%\br It is well known that  And we know weakly quasisymmetric no need to be $\eta$-quasisymmetric, so our result is new even in Euclidean space.
%\er
%\begin{thm}\label{thm1.1}
%Suppose that $X$ and $Y$ are both $c$-quasiconvex complete metric spaces, and $G\varsubsetneq X$ is a domain,
%$G'\varsubsetneq Y$ is a $a$-uniform domain, $f: G\to G'$ and $f^{-1}$ are both $q$-locally $\eta$-quasisymmetric.
%For each subdomain $D$ in $G$, if $D$ is $b$-uniform, then $D'=f(D)$ is $c'$-uniform, where $c'=c'(a,b,c,\varphi)$.
%
%\end{thm}

 Here and in what follows,
the phrase ``the given data  of $X$, $Y$, $G$, $G'$, $D$, $f$ and $f^{-1}$" means the data which depends on the given constants which are the coefficients of quasiconvexity of $X$, $Y$ and $G'$, the coefficient of uniformity of $G$
and  the coefficients of weak quasisymmetry of $f$ and $f^{-1}$.

\br It is worth to mention that in Theorem \ref{thm1.1}, the domain $G'$ is not required to be ``uniform", and only to be ``quasiconvex" (From the definitions
in Section \ref{sec-2}, we easily see that uniformity implies quasiconvexity). Moreover, we see from the example constructed in the paragraph next to Theorem 1.1 in \cite{Xie} that the assumption ``quasiconvexity" in Theorem \ref{thm1.1} is necessary.
\er

 As we have indicated in the second paragraph that
 quasisymmetry and quasiconformality are quantitatively equivalent for homeomorphisms between $\IR^n$. Also, it follows from \cite{Geo} and \cite{TV} that quasiconformality and free quasiconformality are quantitatively equivalent for homeomorphisms between domains in $\IR^n$.
These facts together with \cite[Theorem 1]{HVW} prompt us to conjecture that whether Theorem \ref{thm1.1} still holds if we replace the assumption ``weakly quasisymmetric mappings" (resp. ``$G'$ being quasiconvex") by the one ``freely quasiconformal mappings" (resp. ``$G'$ being uniform"). Under the extra assumption ``$G$  and $G'$ being locally John", we have the following partial answer to this problem.

\begin{thm}\label{thm1.2}
Suppose that $X$ and $Y$ are quasiconvex and complete metric spaces, that
$G \varsubsetneq X$ is a non-point-cut and locally John domain, and that $G'\varsubsetneq Y$ is an uniform and locally John domain. If $f:$ $G \to G'$ is a freely quasiconformal mapping, then for each uniform subdomain $D$ in $G$,
its image $D'=f(D)$ must be uniform, where the coefficient of uniformity of $D'$ depends
only on the given data of $X$, $Y$, $G$, $G'$, $D$ and $f$.
\end{thm}
\noindent Here $G$ is said to be {\it  non-point-cut} if for any $x\in G$, the set $G\backslash\{x\}$ is a subdomain of $G$.

\br When $X=Y=\IR^n$, Theorem \ref{thm1.2} is the subinvariance of uniform domains with respect to quasiconformal mappings in $\IR^n$ (\cite[Remark, P.121]{FHM} and \cite[Theorem 5.6]{Vai2}). Also,
when $X$ and $Y$ are Banach spaces with dimension at least $2$, Theorem \ref{thm1.2} coincides with
\cite[Theorem 1]{HVW}.\er

Also, we get the following subinvariance of uniform domains with respect to quasihyperbolic mappings.

\begin{thm}\label{cor-1}
Suppose that $X$ and $Y$ are quasiconvex and complete metric spaces, that
$G \varsubsetneq X$ is a non-point-cut and locally John domain, and that $G'\varsubsetneq Y$ is an uniform and  locally John domain. If $f:$ $G \to G'$ is a quasihyperbolic mapping, then for each uniform subdomain $D$ in $G$,
its image $D'=f(D)$ must be uniform, where the coefficient of uniformity of $D'$ depends
only on the given data of $X$, $Y$, $G$, $G'$, $D$ and $f$.
\end{thm}

\br When $X$ and $Y$ are Banach spaces with dimension at least $2$, Theorem \ref{cor-1} coincides with Theorem
\Ref{ThmB}, i.e., Theorem 2.44 in \cite{Vai4}.\er

We also conjecture that whether there is the subinvariance of John domains in suitable metric spaces with respect to weakly quasisymmetric mappings, freely quasiconformal mappings etc. See \cite{Heinonen, hlpw} etc for the related discussions in $\IR^n$.

The rest of this paper is organized as follows.  In Section \ref{sec-2}, we recall some definitions and preliminary results. In Section \ref{sec-4}, the proof of Theorem \ref{thm1.1} is given. Section \ref{sec-5} is devoted to the proof of Theorem \ref{thm1.2}, and the proof of Theorem \ref{cor-1} is presented in Section \ref{sec-6}.

\section{Preliminaries}\label{sec-2}

In this section, we give the necessary definitions and auxiliary results, which will be used in the proofs of our main results.

Throughout this paper,
balls and spheres in metric spaces $X$ are written as
$$\mathbb{B}(a,r)=\{x\in X:\,|x-a|<r\},\;\;\mathbb{S}(a,r)=\{x\in X:\,|x-a|=r\}$$
and $$
\mathbb{\overline{B}}(a,r)=\mathbb{B}(a,r)\cup \mathbb{S}(a,r)= \{x\in X:\,|x-a|\leq r\}.$$

For convenience, given
domains $G \subset X,$   $G' \subset Y$, a map $f:G \to G'$
and points $x$, $y$,
$z$, $\ldots$ in  $G$, we always  denote by $x'$, $y'$, $z'$, $\ldots$
the images in $G'$ of $x$, $y$, $z$, $\ldots$ under $f$,
respectively. Also, we assume that $\gamma$
denotes an arc in $G$ and $\gamma'$
 the image in $G'$ of $\gamma$
under $f$.

\subsection{Quasihyperbolic metric, solid arcs and  short arcs}
In this subsection, we start with the definition of quasihyperbolic metric. If $X$ is a connected metric space and $G\varsubsetneq X$ is a non-empty open set,
then it follows from \cite[Remark 2.2]{HL} that the boundary of $G$ satisfies $\partial G\not=\emptyset$.
Suppose $\gamma\subset G$ denotes a rectifiable arc or a path, its {\it quasihyperbolic length} is the number:

$$\ell_{k_G}(\gamma)=\int_{\gamma}\frac{|dz|}{\delta_{G}(z)},
$$ where $\delta_G(z)$ denotes the distance from $z$ to $\partial G$.

%If $\gamma$ is a rectifiable curve with endpoints $x$ and $y$, then
%$$\ell_{k_G}(\gamma)\geq \log\Big(1+\frac{\ell(\gamma)}
%{\min\{\delta_{G}(x), \delta_{G}(y)\}}\Big).
%$$

For each pair of points $x$, $y$ in $G$, the {\it quasihyperbolic distance}
$k_G(x,y)$ between $x$ and $y$ is defined in the following way:
$$k_G(x,y)=\inf\ell_{k_G}(\gamma),
$$
where the infimum is taken over all rectifiable arcs $\gamma$
joining $x$ to $y$ in $G$.

Suppose $X$ is quasiconvex and $G\subsetneq X$. If $\gamma$ is a rectifiable curve in $G$ connecting $x$ and $y$, then (see, e.g., the proof of Theorem $2.7$ in \cite{HL})
$$\ell_{k_G}(\gamma)\geq \log\Big(1+\frac{\ell(\gamma)}
{\min\{\delta_{G}(x), \delta_{G}(y)\}}\Big)
$$
and thus,
$$k_G(x,y)\geq \log\Big(1+\frac{|x-y|}
{\min\{\delta_{G}(x), \delta_{G}(y)\}}\Big).
$$

Gehring and Palka \cite{GP} introduced the quasihyperbolic metric of
a domain in $\IR^n$. For the basic properties of this metric we refer to \cite{Geo}. Recall that a curve $\gamma$ from $x$ to
$y$ is a {\it quasihyperbolic geodesic} if
$\ell_{k_G}(\gamma)=k_G(x,y)$. Each subcurve of a quasihyperbolic
geodesic is obviously a quasihyperbolic geodesic. It is known that a
quasihyperbolic geodesic between any two points in a Banach space $X$ exists if the
dimension of $X$ is finite, see \cite[Lemma 1]{Geo}. This is not
true in arbitrary metric spaces (cf. \cite[Example 2.9]{Vai6-0}).

%%%%%%%%%%%%%%%%%
Let us recall a result which is useful for the discussions later on.

\begin{Lem}\label{ll-11}(\cite[Lemma 2.4]{HWZ}) Let $X$ be a $c$-quasiconvex metric space and let $G\subsetneq X$ be a domain. Suppose that $x$, $y\in G$ and either $|x-y|\leq \frac{1}{3c}\delta_G(x)$ or $k_G(x,y)\leq 1$. Then
\be\label{vvm-2} \frac{1}{2}\frac{|x-y|}{\delta_G(x)}< k_G(x,y) < 3c\frac{|x-y|}{\delta_G(x)}.\ee
\end{Lem}

\noindent Here, we say that $X$ is {\it c-quasiconvex} $(c\geq 1)$ if each pair of points $x$, $y\in X$ can be joined by an arc $\gamma$ with length ${\ell}(\gamma)\leq c|x-y|$.

\bdefe \label{def1.4}
 Suppose $\gamma$ is an arc in a domain $G\varsubsetneq X$ and $X$ is a rectifiably connected metric space. The arc may be closed, open or half open. Let $\overline{x}=(x_0,$ $\ldots,$ $x_n)$,
$n\geq 1$, be a finite sequence of successive points of $\gamma$.
For $h\geq 0$, we say that $\overline{x}$ is {\it $h$-coarse} if
$k_G(x_{j-1}, x_j)\geq h$ for all $1\leq j\leq n$. Let $\Phi_k(\gamma,h)$
denote the family of all $h$-coarse sequences of $\gamma$. Set

$$s_k(\overline{x})=\sum^{n}_{j=1}k_G(x_{j-1}, x_j)$$ and
$$\ell_{k_G}(\gamma, h)=\sup \{s_k(\overline{x}): \overline{x}\in \Phi_k(\gamma,h)\}$$
with the agreement that $\ell_{k_G}(\gamma, h)=0$ if
$\Phi_k(\gamma,h)=\emptyset$. Then the number $\ell_{k_G}(\gamma, h)$ is the
{\it $h$-coarse quasihyperbolic length} of $\gamma$.  \edefe

If $X$ is $c$-quasiconvex, then $\ell_{k_G}(\gamma, 0)=\ell_{k_G}(\gamma)$ (see, e.g., \cite[Proposition A.7 and Remark A.13]{BHK} and \cite[Lemma 2.5]{HWZ} ).

\bdefe \label{def1.5} Let $G$ be a proper domain in a rectifiably connected metric space $X$. An arc $\gamma\subset G$
is {\it $(\nu, h)$-solid} with $\nu\geq 1$ and $h\geq 0$ if
$$\ell_{k_G}(\gamma[x,y], h)\leq \nu\;k_G(x,y)$$ for all $x$, $y\in \gamma$.

An arc $\gamma\subset G$ with endpoints $x$ and $y$ is said to be $\varepsilon$-short ($\varepsilon\geq 0$) if $$\ell_{k_G}(\gamma)\leq k_G(x,y)+\varepsilon.$$

Obviously, by the definition of $k_G$, we know that for every $\varepsilon>0$, there exists an arc $\gamma \subset G$ such that $\gamma$ is $\varepsilon$-short, and it is easy to see that  every subarc of an $\varepsilon$-short arc is also $\varepsilon$-short.\edefe

\br
For any pair of points $x$ and $y$ in a proper domain $G$ of Banach space $E$, if the dimension of $E$ is finite, then there exists a quasihyperbolic geodesic in $G$ connecting $x$ and $y$ (see \cite[Lemma 1]{Geo}).  But if the dimension of $E$ is infinite, this property is no longer valid (see, e.g., \cite[Example 2.9]{Vai6-0}). In order to overcome this shortcoming in Banach spaces, V\"ais\"al\"a proved the existence of neargeodesics or quasigeodesics (see \cite{Vai6}), and  every quasihyperbolic geodesic is a quasigeodesic. See also \cite{RT}. In metric spaces, we do not know if this existence property is true or not. However, this existence property plays a very important role in the related discussions.
In order to overcome this disadvantage, in this paper, we will exploit the substitution of ``quasigeodesics" replaced by ``short arcs". The class of short arcs has been introduced when V\"ais\"al\"a studied properties of Gromov hyperbolic spaces \cite{Vai9} (see also \cite{BH,Herron}), and as we see that the existence of such class of arcs is obvious in metric spaces.
Although, there is no implication between the class of quasigeodesics and the one of short arcs, we will prove that, under certain geometric assumptions, every short arc is a double cone arc.
\er

%%%%%%%%%%%%%%%%%

By a slight modification of the method used in the proof of \cite[Lemma 6.21]{Vai6}, we get the following  result.

\begin{lem}\label{ll-14} Suppose that $X$ is a $c$-quasiconvex metric space and that $G\varsubsetneq X$ is a domain, and that $\gamma$ is a  $(\nu,h)$-solid arc in $G$ with endpoints $x$, $y$ such that $\min\{\delta_G(x),\delta_G(y)\}=r\geq 3c|x-y|$. Then there is a constant $\mu_1=\mu_1(c,\nu)$ such that $$\diam(\gamma)\leq \max\{\mu_1|x-y|, 2r(e^h-1)\},$$ where ``$\diam$" means ``{\rm diameter}".
\end{lem}

\bpf Without loss of generality, we assume that $\delta_G(y)\geq \delta_G(x)=r$. Denoting $t=|x-y|$ and applying Lemma \Ref{ll-11}, we get $$k_G(x,y)\leq 3ct/r.$$
Let $u\in \gamma$. To prove this lemma, it suffices to show that there exists a constant $\mu_1=\mu_1(c,\nu)$ such that
\be\label{neq-eq-1}|u-x|\leq \max\big\{\frac{\mu_1}{2}|x-y|,r(e^h-1)\big\}.\ee

To this end, we consider two cases. The first case is:  $k_G(u,x)\leq h$. Under this assumption, it is easy to see that
\be\label{dw-2}|u-x|\leq (e^{k_G(u,x)}-1)\delta_G(x)\leq r(e^h-1).\ee

For the remaining case: $k_G(u,x)> h$, we choose a sequence of successive points of $\gamma$: $x=x_0$, $\ldots$, $x_{n}=u$ such that $$ k_G(x_{j-1},x_j)=h \;\;\; {\rm for}\;\;\; j\in\{1,\ldots, n-1\}$$ and $$0< k_G(x_{n-1},x_n)\leq h.$$ Then $n\geq 2$ and
$$(n-1)h\leq \sum_{i=1}^{n-1}k_G(x_{j-1},x_j)\leq \ell_{k_G}(\gamma,h)\leq \nu k_G(x,y)\leq 3c\nu t/r,$$
which shows that $$k_G(x,u)\leq \sum_{i=1}^{n}k_G(x_{j-1},x_j)\leq nh\leq 6c\nu t/r.$$
Let $s=\frac{t}{r}$. Then $s\leq \frac{1}{3c}$ and
$$\frac{|u-x|}{t}\leq \frac{e^{6c\nu s}-1}{s}.$$
Obviously,
the function $g(s)=\frac{1}{s}\big(e^{6c\nu s}-1\big)$ is increasing in
$(0,\frac{1}{3c}]$ and $\lim_{s\to 0}{\frac{e^{6c\nu s}-1}{s}}=6c\nu $. Letting $$\mu_1=6c(e^{2\nu}-1)$$ gives
\be\label{dw-1} |u-x|\leq \frac{1}{2}\mu_1t.\ee
It follows from \eqref{dw-2} and \eqref{dw-1} that \eqref{neq-eq-1} holds, and hence the proof of the lemma is complete.
\epf

\begin{lem} \label{mon-4}  Suppose that $X$ is a $c$-quasiconvex metric space and $G\subsetneq X$ is a domain. Suppose, further, that for $x$,
$y\in G$, \begin{enumerate}
\item
$\gamma$ is an $\varepsilon$-short arc in $G$ connecting $x$ and $y$ with $0<\varepsilon\leq \frac{1}{2}k_{G}(x,y)$, and
\item
$|x-y|\leq \frac{1}{3c} \min\{\delta_{G}(x), \delta_{G}(y)\}$. \end{enumerate}
Then $$\ell(\gamma)\leq \frac{9}{2}ce^{\frac{3}{2}}|x-y|.$$\end{lem}
\bpf Without loss of generality, we assume that $\min\{\delta_{G}(x), \delta_{G}(y)\}=\delta_{G}(x)$. It follows from Lemma \Ref{ll-11} that
$$\log\left(1+\frac{\ell(\gamma)}{\delta_{G}(x)}\right)\leq \ell_{k_{G}}(\gamma)\leq k_{G}(x,y)+\varepsilon \leq  \frac{3}{2}k_{G}(x,y)\leq \frac{9c}{2}\frac{|x-y|}{\delta_{G}(x)}\leq \frac{3}{2}.
$$ Hence,
$$\frac{\ell(\gamma)}{\delta_{G}(x)}\leq e^{\frac{3}{2}}-1,$$
which leads to
$$\frac{\ell(\gamma)}{\delta_{G}(x)}\leq e^{\frac{3}{2}}\log\left(1+\frac{\ell(\gamma)}{\delta_{G}(x)}\right)\leq \frac{9}{2}ce^{\frac{3}{2}}\frac{|x-y|}{\delta_{G}(x)}.$$ Therefore, $$\ell(\gamma)\leq \frac{9}{2}ce^{\frac{3}{2}}|x-y|,$$ as required.\epf

\subsection{Uniform domains, John domains and  locally John domains}
In 1961, John \cite{John} introduced the twisted interior cone condition with his work on elasticity, and these domains where first called John domains by Martio and Sarvas in \cite{MS}. In the same paper, Martio and Sarvas  also introduced another class of domains which are the uniform domains. The main motivation for studying these domains was in showing global injectivity properties for locally injective mappings. Since then, many other characterizations of uniform and John domains have been established, see
\cite{FW,  Geo, Martio-80, Vai6, Vai4, Vai8}, and the importance of these classes of domains in the function theory is well
documented (see e.g. \cite{FW, Vai2}). Moreover, John and uniform domains in
$\mathbb{R}^n$ enjoy numerous geometric and function theoretic
properties that are useful in other many fields of modern mathematical analysis as well (see e.g.
\cite{Alv, Bea, Geo, Has, Ml, Jo80, Yli, Vai2}, and references therein).

We recall the definitions of uniform domains and John domains following closely the notation
and terminology of \cite{TV, Vai2, Vai, Vai6-0, Vai6}  and \cite{Martio-80}.

\bdefe \label{def1.3} A domain $G$ in $X$ is called $b$-{\it
uniform} provided there exists a constant $b$
with the property that each pair of points $x$, $y$ in $G$ can
be joined by a rectifiable arc $\gamma$ in $G$ satisfying
 \bee
\item $\min\{\ell(\gamma[x,z]),\ell(\gamma[z,y])\}\leq b\,\delta_{G}(z)$ for all $z\in \gamma$, and
\item $\ell(\gamma)\leq b\,|x-y|$,
\eee
\noindent where $\ell(\gamma)$ denotes the length of $\gamma$,
$\gamma[x,z]$ the part of $\gamma$ between $x$ and $z$.
At this time, $\gamma$ is said to be a {\it double $b$-cone arc}. Condition $(1)$ (resp. $(2)$) is called the {\it cigar condition} (resp. {\it  turning condition}).

If the condition $(1)$ is satisfied, not necessarily $(2)$, then $G$ is said to be a {\it $b$-John domain}. At this time, the arc $\gamma$ is called a {\it $b$-cone arc}.
\edefe
\bdefe\label{wed-1}
 A domain $G$ in $X$ is said to be a {\it locally $a$-John domain} if there exists a constant $a$ such that all metric balls $\mathbb{B}(x, r)$ are $a$-John domains,
where $x\in G$ and $0<r\leq \delta_G(x)$. \edefe

We note that all domains which satisfies the strong geodesic condition (see \cite{BKL} for the definition) in the abstract setting of homogeneous spaces are locally John domains \cite[Corollary 3.2]{BKL}. In particular, all domains in Carnot-Carath\'eodory metric spaces \cite{BKL} and Banach spaces \cite{Vai6} are locally John domains.

We remark that for $x\in G$, it is possible that the metric ball $\mathbb{B}(x,r)$ is not contained in $G$ for $0<r\leq \delta_G(x)$. But, in \cite{HWZ},
the authors proved the following.

\begin{Lem}\label{Lem-J} Suppose $X$ is a $c$-quasiconvex metric space and $G\subsetneq X$ is a domain. Then for any rectifiably connected open set
$D\subset \mathbb{B}(x,r)$ with $x\in D\cap G$, if $r\leq \delta_G(x)$, then $D \subset G.$
\end{Lem}

It follows from Lemma \Ref{Lem-J} that the following result is obvious.

\begin{lem}\label{lem-LHWZ1}
Suppose $X$ is a $c$-quasiconvex metric space and $G\subsetneq X$ is  a locally John domain.
For any $x\in G$, if $0<r\leq \delta_G(x)$, then $\mathbb{B}(x,r)\subset G$.
\end{lem}

Further, for locally John domains, we have the following estimate for the quasihyperbolic metric.

\begin{lem}\label{LHWZ-lem1} Suppose $G\varsubsetneq X$ is a locally $a$-John domain. Then for $x,$ $y\in G$ with $|x-y|=t \delta_G(x)$, where $0<t<1$, we have $$k_G(x,y) \leq 2a\frac{3+t}{1-t}.$$
\end{lem}
\bpf
Let $x,$ $y\in G$ with $|x-y|=t \delta_G(x)$. Since $G$ is locally $a$-John, we know from Definition \ref{wed-1} and Lemma \ref{LHWZ-lem1} that $\mathbb{B}(x, \frac{1+t}{2}\delta_G(x))\subset G$ is an $a$-John domain. Then  there is a curve $\gamma\in \mathbb{B}(x,  \frac{1+t}{2}\delta_G(x))$ such that for any $w\in \gamma$, $$\min\{\ell(\gamma[y, w]),\ell(\gamma[x,w])\}\leq a \delta_G(w)$$ and $$\delta_G(w)\geq \delta_G(x)-|x-w|\geq \frac{1}{2}(1-t)\delta_G(x).$$
Let $x_0\in \gamma$ be the point bisecting the arc length of $\gamma$. Then $$\ell(\gamma)\leq 2a\delta_G(x_0)\leq 2a(\delta_G(x)+|x-x_0|)\leq a(3+t)\delta_G(x),$$ which shows that  $$k_G(x,y)\leq \int_{\gamma}\frac{|dw|}{\delta_G(w)}\leq \frac{2\ell(\gamma)}{(1-t)\delta_G(x)}
\leq 2a\frac{3+t}{1-t},$$ as required.
\epf

Let us recall the following useful property of uniform domains.

\begin{Lem}\label{BHK-lem}$($\cite[Lemma 3.12]{BHK}$)$ Suppose $G\subsetneq X$ is a $b$-uniform domain in a rectifiable connected metric space $X$. Then for any $x, y\in G$, we have $$k_G(x,y)\leq 4b^2\log\left(1+\frac{|x-y|}{\min\{\delta_G(x),\delta_G(y)\}}\right).$$
\end{Lem}

The following are the analogues of Lemmas $6.10$ and $6.11$ in \cite{Vai6} in the setting of metric spaces. The proofs are similar.

\begin{Lem}\label{ll-12} Suppose that $G\subsetneq X$ is a $b$-uniform domain in a rectifiable connected metric space $X$, and that $\gamma$ is an arc in $\{x\in G: \delta_G(x)\leq r\}$. If $\gamma$ is $(\nu,h)$-solid, then  $$\diam(\gamma)\leq M_1 r,$$ where $M_1=M_1(b, \nu,h)$.
\end{Lem}

%%%%%%%%%%%%%%%%%

\begin{Lem}\label{ll-13} For all $b\geq 1$, $\nu\geq 1$ and $h\geq 0$, there are constants $0<q_0=q_0(b,\nu,h)<1$ and $M_2=M_2(b,\nu,h)\geq 1$
with the following property:
Suppose that $G$ is a $b$-uniform domain and $\gamma$ is a $(\nu,h)$-solid arc starting at $x_0\in G$. If $\gamma$ contains a point $u$ with $\delta_G(u)\leq q_0\delta_G(x_0)$, then $$\diam(\gamma_{u})\leq M_2\delta_G(u),$$ where $\gamma_{u}=\gamma\setminus \gamma[x_0,u)$.
\end{Lem}

%%%%%%%%%%%%%%%%%

Now, we are ready to prove an analogue of Lemma \ref{ll-14} for uniform domains.

\begin{lem}\label{ll-15} Suppose that $X$ is a $c$-quasiconvex metric space and that $G\varsubsetneq X$ is a $b$-uniform domain, and that $\gamma$ is a $(\nu,h)$-solid arc in $G$ with endpoints $x$, $y$. Let $\delta_G(x_0)=\max_{p\in \gamma}\delta_G(p)$. Then there exist constants
$\mu_2=\mu_2( b, \nu, h)\geq 1$ and $\mu_3=\mu_3( b, c, \nu, h)\geq 1$ such that
\begin{enumerate} \item\label{ma-0-1}
$\diam(\gamma[x,u])\leq \mu_2 \delta_G(u)$ for $u\in \gamma[x,x_0],$ and $\diam(\gamma[y,v])\leq \mu_2 \delta_G(v)$ for $v\in \gamma[y, x_0]$;
\item\label{ma-0-2} $\diam(\gamma)\leq \max\big\{\mu_3 |x-y|, 2(e^h-1)\min\{\delta_G(x),\delta_G(y)\}\big\}.$
\end{enumerate}
 \end{lem}

\bpf We first prove \eqref{ma-0-1}. Obviously, it suffices to prove the first inequality in \eqref{ma-0-1} because the proof for the second one is similar. Let $$\mu_2=\max\Big\{\frac{M_1}{q_0}, M_2\Big\},$$ where $M_1=M_1(b,\nu,h)$ is the constant from Lemma
\Ref{ll-12}, $q_0=q_0(b,\nu,h)$ and $M_2=M_2(b,\nu,h)$ are the constants from Lemma \Ref{ll-13}.

For $u\in \gamma[x,x_0]$, we divide the proof into two cases. If $\delta_G(u)\leq q_0\delta_G(x_0)$, then Lemma \Ref{ll-13} leads to
\be\label{ma-1} \diam(\gamma[x,u])\leq M_2\delta_G(u).\ee If  $\delta_G(u)> q_0\delta_G(x_0)$, then applying Lemma \Ref{ll-12} with the substitution
$r$ replaced by $\delta_G(x_0)$ and $\gamma$ replaced by $\gamma[x,u]$, we easily get
\be\label{ma-2}\diam(\gamma[x,u])\leq M_1\delta_G(x_0)< \frac{M_1}{q_{0}}\delta_G(u).\ee
It follows from \eqref{ma-1} and \eqref{ma-2} that the first assertion in \eqref{ma-0-1} holds, and thus the proof of \eqref{ma-0-1} is complete.

To prove \eqref{ma-0-2}, without loss of generality, we assume that $$\min\{\delta_G(x),\; \delta_G(y)\}= \delta_G(x)\;\;\; {\rm and}\;\;\; \diam(\gamma)> |x-y|.$$
Let
$$\mu_3=\frac{3}{4}\big[1+2(1+6c)(e^{h+4b^2\nu\log(1+4\mu_2)}-1)\big].$$

If $\delta_G(x)\geq 3c|x-y|$, then $(2)$ follows from Lemma \ref{ll-14} since the constant $\mu_1$ in Lemma \ref{ll-14} satisfies $\mu_1< \mu_3$.
Hence, in the following, we assume that $$\delta_G(x)<3c|x-y|.$$ Let $x_1\in \gamma$ (resp. $y_1\in \gamma$) be the first point in $\gamma$ from $x$ to $y$ (resp. from $y$ to $x$) such that
\beq\label{eq-new-1}\diam(\gamma[x,x_1])=\frac{1}{2}|x-y|\;\,(\mbox{resp.}\; \diam(\gamma[y,y_1])=\frac{1}{2}|x-y|).\eeq

%\begin{figure}[htbp]
%\begin{center}
%\input{figure01.pspdftex}
%\caption{The points $w_1$, $x_1,$ $y_1$ and $w_2$ in $\gamma$.}
%
%\label{subinv-fig01}
%\end{center}
%\end{figure}

%{\color{red}Figure $1$: The points $w_1$, $x_1$, $y_1$ and $w_2$ in $\gamma$.}

\noindent Then we have $$\diam(\gamma[y,x_1])\geq |y-x_1|\geq |y-x|-|x-x_1|\geq\frac{1}{2}|y-x|=\diam(\gamma[x,x_1]),$$ and similarly, we get $$\diam(\gamma[x,y_1])>\diam(\gamma[y,y_1]).$$
Thus, it follows from  \eqref{ma-0-1} that $$\frac{1}{2}|x-y|=\diam(\gamma[x,x_1])=\diam(\gamma[y,y_1])\leq  \mu_2
\min\{\delta_G(x_1), \delta_G(y_1)\}.$$ Also,
$$|x_1-y_1|\leq |x_1-x|+|x-y|+|y-y_1|\leq 2|x-y|.$$ Then Lemma \Ref{BHK-lem} implies
$$k_G(x_1,y_1)\leq 4b^2\log\left(1+\frac{|x_1-y_1|}{\min\{\delta_G(x_1),\delta_G(y_1)\}}\right)\leq 4b^2\log(1+4\mu_2).$$
Since $\gamma$ is a $(\nu,h)$-solid arc, for any $u_1,$ $u_2\in \gamma[x_1,y_1]$, we have
\begin{eqnarray*}k_G(u_1,u_2)&\leq&   \max\{h,\; \ell_{k_G}(\gamma[x_1,y_1],h)\}\leq h+\nu k_G(x_1,y_1)\\&\leq& h+4b^2\nu\log(1+4\mu_2),\end{eqnarray*}
and so, for all $z\in \gamma[x_1,y_1]$,
\begin{eqnarray}\label{eq-new-2}|z-x_1|&\leq& (e^{k_G(z,x_1)}-1)\delta_G(x_1)\\\nonumber&\leq&   (e^{h+4b^2\nu\log(1+4\mu_2)}-1)(\delta_G(x)+|x-x_1|)\\\nonumber&\leq& \frac{1}{2}(1+6c)(e^{h+4b^2\nu\log(1+4\mu_2)}-1)|x-y|.\end{eqnarray}
Let $w_1,$ $w_2\in \gamma$ be points such that \be \label{sat-1}|w_1-w_2|\geq \frac{2}{3}\diam(\gamma).\ee
Then we get
\bcl\label{sat-2}
$|w_1-w_2|\leq \frac{2}{3}\mu_3|x-y|$.
\ecl
Since \eqref{eq-new-1} guarantees that neither $\gamma[x,x_1]$ nor $\gamma[y, y_1]$ contains the set $\{w_1,w_2\}$, we see that,
to prove this claim, according to the positions of $w_1$ and $w_2$ in $\gamma$, we need to consider the following four possibilities.\begin{enumerate}
\item
 $w_1\in \gamma[x,x_1]$ and $w_2\in \gamma[y,y_1]$.
 Obviously, by \eqref{eq-new-1}, we have  $$|w_1-w_2|\leq |w_1-x|+|x-y|+|y-w_2|\leq 2|x-y|.$$

\item
 $w_1\in \gamma[x,x_1]$ and $w_2\in \gamma[x_1,y_1]$. Then \eqref{eq-new-1} and \eqref{eq-new-2} show that $$|w_1-w_2|\leq |w_1-x_1|+|x_1-w_2|\leq \frac{1}{2}\big[1+(1+6c)(e^{h+4b^2\nu\log(1+4\mu_2)}-1)\big]|x-y|.$$

\item
$w_1,$ $w_2\in \gamma[x_1,y_1]$. Then \eqref{eq-new-2} implies $$|w_1-w_2|\leq |w_1-x_1|+|x_1-w_2|\leq (1+6c)(e^{h+4b^2\nu\log(1+4\mu_2)}-1)|x-y|.$$

\item
$w_1\in \gamma[x_1,y_1]$ and $w_2\in \gamma[y_1,y]$. Again, we infer  from \eqref{eq-new-1} and \eqref{eq-new-2} that $$ |w_1-w_2|\leq |w_1-x_1|+|x_1-y_1|+|y_1-w_2|\leq \frac{1}{2}\big[1+2(1+6c)(e^{h+4b^2\nu\log(1+4\mu_2)}-1)\big]|x-y|.$$
\end{enumerate}
The claim is proved.\medskip

Now, we are ready to finish the proof. It follows from \eqref{sat-1} and Claim \ref{sat-2} that
$$\diam(\gamma)\leq \frac{3}{2}|w_1-w_2|\leq \mu_3|x-y|,$$
which implies that \eqref{ma-0-2} also holds in this case. Hence, the proof of the lemma is complete.
\epf

\subsection{Quasihyperbolic mappings, coarsely quasihyperbolic mappings and relative homeomorphisms}

\bdefe Let $G\varsubsetneq X$ and $G'\varsubsetneq Y$ be two domains. We say
that a homeomorphism $f: G\to G'$ is \begin{enumerate}

\item
 {\it $C$-coarsely $M$-quasihyperbolic}, or briefly
$(M,C)$-CQH, if there are constants $M\geq 1$ and $C\geq 0$ such that for all $x$, $y\in G$,
$$\frac{k_G(x,y)-C}{M}\leq k_{G'}(f(x),f(y))\leq M\;k_G(x,y)+C.$$

\item {\it $M$-quasihyperbolic},
or briefly  {\it $M$-QH}, if $f$ is $(M,0)$-CQH.
\item
{\it fully $C$-coarsely $M$-quasihyperbolic}  if there are constants $M\geq 1$ and $C\geq 0$ such that
$f$ is $C$-coarsely $M$-quasihyperbolic in every subdomain of $G$.

\end{enumerate}\edefe

Under coarsely quasihyperbolic mappings, we have the following useful relationship between short arcs and solid arcs.

\begin{lem}\label{ll-001} Suppose that $X$ and $Y$ are rectifiable connected metric spaces, and that $G\varsubsetneq X$ and $G'\varsubsetneq Y$ are domains. If $f:\;G\to G'$ is $(M,C)$-CQH, and $\gamma$ is an $\varepsilon$-short arc in $G$ with $0<\varepsilon\leq 1$, then there are constants $\nu=\nu(C,M)$ and $h=h(C, M)$ such that
 the image $\gamma'$ of $\gamma$ under $f$ is $(\nu,h)$-solid in $G'$.
\end{lem}
\bpf
Let$$h=(2M+1)C+2M\;\; \mbox{and}\;\; \nu=\frac{4(C+1)M(M+1)}{2C+1}.$$
Obviously, we only need to verify that for $x$ and $y\in \gamma$,
\be\label{new-eq-3}\ell_{k_{G'}}(\gamma'[x',y'],h)\leq\nu k_{G'}(x',y').\ee We prove this by considering two cases.
The first case is: $k_G(x,y)<2C+1$. Then for $z_1$, $z_2\in\gamma[x, y]$, we have
$$k_{G'}(z'_1,z'_2)\leq Mk_G(z_1,z_2)+C\leq M(k_G(x,y)+\varepsilon)+C<(2M+1)C+2M=h,$$
and so
\be\label{ma-3}\ell_{k_{G'}}(\gamma'[x',y'],h)=0.\ee

Now, we consider the other case: $k_G(x,y)\geq 2C+1$. Then $$k_{G'}(x',y')\geq \frac{1}{M}(k_G(x,y)-C)> \frac{1}{2M}k_G(x,y).$$ With the aid of \cite[Theorems 4.3 and 4.9]{Vai6}, we have
\beq\label{ma-4}
\ell_{k_{G'}}(\gamma'[x',y'],h) &\leq&
\ell_{k_{G'}}(\gamma'[x',y'],(M+1)C) \leq
 (M+1)\ell_{k_G}(\gamma[x,y])\\ \nonumber &\leq&(M+1)(k_G(x,y)+\varepsilon)\\ \nonumber &\leq& \frac{2(C+1)(M+1)}{2C+1}k_{G}(x,y) \\ \nonumber &\leq& \frac{4(C+1)M(M+1)}{2C+1}k_{G'}(x',y').\eeq
It follows from \eqref{ma-3} and \eqref{ma-4} that \eqref{new-eq-3} holds.\epf

%%%%%%%%%%%%%%%%%

The following two lemmas are useful in the proof of Theorem \ref{thm1.1}.

\begin{lem}\label{ll-000} Suppose that $X$ and $Y$ are both $c$-quasiconvex and complete metric spaces, and that $G\varsubsetneq X$ and $G'\varsubsetneq Y$ are domains. If both $f:$ $G\to G'$ and $f^{-1}:$ $G'\to G$ are weakly $H$-quasisymmetric, then
\begin{enumerate}
\item \label{sat-3}
$f$ is $\varphi$-FQC, where $\varphi=\varphi_{c,H}$ which means that the function $\varphi$ depends only on $c$ and $H$;

\item \label{sat-4}
 $f$ is fully $(M,C)$-CQH, where $M=M(c,H)\geq 1$ and $C=C(c,H)\geq 0$ are constants.
 \end{enumerate}
\end{lem}
\bpf By \cite[Theorem 1.6]{HL}, we know that for every subdomain $D\subset G$, both $f:$ $D\to D'$ and $f^{-1}:$ $D'\to D$ are $\varphi$-semisolid with $\varphi=\varphi_{c,H}$, and so, $f$ is $\varphi$-FQC. Hence \eqref{sat-3} holds.  Meanwhile,  \cite[Theorem 1]{ HWZ} implies that \eqref{sat-3} and \eqref{sat-4} are equivalent, and thus, \eqref{sat-4} also holds.
\epf

%%%%%%%%%%%%%%%%%
\begin{Lem}\label{lem-ll-0}$($\cite[Lemma 6.5]{Vai8}$)$ Suppose that $X$ is $c$-quasiconvex, and that $f:$ $X\to Y$ is weakly $H$-qasisymmetric. If $x,$ $y,$ $z$ are distinct points in $X$ with $|y-x|\leq t|z-x|$, then $$|y'-x'|\leq \theta(t)|z'-x'|,$$ where the function $\theta(t)=\theta_{c, H}(t)$ is increasing in $t$.
\end{Lem}

%%%%%%%%%%%%%%%%%
\bdefe  Let $G\varsubsetneq X$ and $G'\varsubsetneq Y$ be two domains. We say
that a homeomorphism $f: G\to G'$ is  \begin{enumerate}
\item
$(\theta, t_{0})$-relative if there is a constant $t_0\in (0, 1]$ and a homeomorphism $\theta: [0, t_0)\to [0, \infty)$ such that $$\frac{|x'-y'|}{\delta_{G'}(x')} \leq\theta\Big(\frac{|x-y|}{\delta_G(x)}\Big)$$
whenever $x$, $y\in G$ and $|x-y|<t_{0}\delta_G(x)$; In particular, if $t_0=1$, then $f$ is called to be $\theta$-relative;
\item
fully $(\theta, t_{0})$-relative (resp. fully $\theta$-relative) if $f$ is $(\theta, t_{0})$-relative (resp. $\theta$-relative) in every subdomain of $G$.
\end{enumerate}\edefe

%%%%%%%%%%%%%%%%%

\section{The proof of Theorem \ref{thm1.1}}\label{sec-4}

In this section, we always assume that $X$ and $Y$ are $c$-quasiconvex and complete metric spaces, and that $G\varsubsetneq X$ and $G'\varsubsetneq Y$ are domains. Furthermore, we suppose that both $f:$ $G\to G'$ and $f^{-1}:$ $G'\to G$ are weakly $H$-quasisymmetric, $G'$ is $c_1$-quasiconvex and $D\subset G$ is $b$-uniform.

Under these assumptions, it follows from Lemma \ref{ll-000} that $f$ is $(M, C)$-$CQH$ with $M=M(c, H)\geq 1$ and $C=C(c, H)\geq 0$.

We are going to show the uniformity of $D'=f(D)$. For this, we let $x'$, $y'\in D'=f(D)\subset G'$, and $\gamma'$ be an $\varepsilon$-short arc in $D'$ joining $x'$ and $y'$ with
$$0<\varepsilon<\min\big\{1,\frac{1}{2}k_{D'}(x',y')\big\}.$$
Then by Lemmas \ref{ll-001} and \ref{ll-000}\eqref{sat-4}, the preimage $\gamma$ of $\gamma'$ is a $(\nu,h)$-solid arc in $D$ with $\nu=\nu(c, H)$ and $h=h(c, H)$. Let $w_0\in\gamma$ be such that
\be\label{wes-1}
\delta_D(w_0)=\max_{p\in \gamma}\delta_D(p).
\ee

%\begin{figure}[htbp]
%\begin{center}
%\input{figure03.pspdftex}
%\caption{The arc $\gamma$ and the related points.}
%
%\label{subinv-fig03}
%\end{center}
%\end{figure}

\noindent Then by Lemma \ref{ll-15}, there is a constant $\mu=\mu(b,\nu,h)$ such that for each $u\in\gamma[x,w_0]$ and for all $z\in\gamma[u,w_0]$, \be\label{new-eq-4}|u-z|\leq \diam (\gamma[u, z])
\leq\mu\delta_D(z),\ee
and for each $v\in\gamma[y,w_0]$ and for all $z\in\gamma[v,w_0]$, \be\nonumber |v-z|\leq \diam (\gamma[v, z])\leq\mu\delta_D(z).\ee

In the following, we show that $\gamma'$ is a double cone arc in $D'$. Precisely, we shall prove
that there exist constants $A\geq 1$ and $B\geq 1$ such that for every $z'\in\gamma'[x',y']$,
\be\label{main-eq-1}\min\{\ell(\gamma'[x',z']),\ell(\gamma'[z',y'])\}\leq A\delta_{D'}(z')\ee
and
\be\label{main-eq-2}\ell(\gamma')\leq B|x'-y'|.\ee

The verification of \eqref{main-eq-1} and \eqref{main-eq-2} is given in the following two subsections.

\subsection{The proof of \eqref{main-eq-1}}

Let
$$A=2e^{8b^2A_1(C+1)M}\;\; {\rm and}\;\; A_1=2e^{M+C}(1+\mu)\theta''\Big(6c\theta'(\mu)e^{4b^2M+C}\Big),$$ where the functions $\theta'=\theta'_{b,H}$ and $\theta''=\theta''_{c_1,H}$ are from Lemma \Ref{lem-ll-0}.  Obviously, we only need to get the following estimate:
For all $z'\in\gamma'[x',w'_0]$ (resp. $z'\in\gamma'[y',w'_0]$),
\be\label{main-eq-3}
\ell(\gamma'[x',z'])\leq A\delta_{D'}(z')\; \; ({\rm resp.}\; \ell(\gamma'[y',z'])\leq A\delta_{D'}(z')).
\ee

It suffices to prove the case $z'\in\gamma'[x',w'_0]$ since the proof of the case $z'\in\gamma'[y',w'_0]$ is similar.
Suppose on the contrary that there exists some point $x'_0\in \gamma'[x',w'_0]$ such that
\be\label{main-eq-4}
\ell(\gamma'[x',x_0'])>A\delta_{D'}(x_0').
\ee Then we choose $x'_1\in\gamma'[x',w'_0]$ be the first point from $x'$ to $w_0'$ such that
\be\label{sat-5}
\ell(\gamma'[x',x'_1])=A\delta_{D'}(x'_1).
\ee

%\begin{figure}[htbp]
%\begin{center}
%\input{figure02.pspdftex}
%\caption{The arc $\gamma'$ and the related points.}
%
%\label{subinv-fig02}
%\end{center}
%\end{figure}

%{\color{red} Figure $2$: The arc $\gamma'$ and the related points.}

Let $x_2\in  D$ be such that
$$|x_1-x_2|=\frac{1}{2}\delta_{D}(x_1).$$

%{\color{red} Figure $3$: The arc $\gamma$ and the related points.}

\noindent Then we have
\bcl\label{eq-lwz4}
$|x_1'-x_2'|< e^{4b^2M+C}\delta_{D'}(x_1').$
\ecl
Obviously,
$$\delta_{D}(x_2)\geq \delta_{D}(x_1)-|x_1-x_2|=|x_1-x_2|,$$
and so, Lemma \Ref{BHK-lem} implies
\begin{eqnarray*}\log\left(1+\frac{|x_1'-x_2'|}{\delta_{D'}(x_1')}\right)&\leq& k_{D'}(x_1',x_2')\leq Mk_D(x_1,x_2)+C\\&\leq& 4b^2M\log\left(1+\frac{|x_1-x_2|}{\min\{\delta_D(x_1),\delta_D(x_2)\}}\right)+C\\&<& 4b^2M+C,
\end{eqnarray*}
whence
$$
|x_1'-x_2'|< e^{4b^2M+C}\delta_{D'}(x_1'),
$$ which shows that the claim holds.\medskip

Let $x'_3\in\gamma'[x',x'_1]$ be such that
\be\label{sun-2}\ell(\gamma'[x',x'_3])=\frac{1}{2}\ell(\gamma'[x',x'_1]),\ee
and then, we get an estimate on $|x_1'-x_2'|$ in terms of $d_{D'}(x_3')$ as stated in the following claim.
\bcl\label{eq-lwz1}
$|x_1'-x_2'|< 2e^{4b^2M+C}\delta_{D'}(x_3').$
\ecl
It follows from \eqref{sat-5} and \eqref{sun-2} that
$$
\delta_{D'}(x_1')< 2 \delta_{D'}(x_3'),
$$ since the choice of $x_1'$ implies $\ell(\gamma'[x',x'_3])<Ad_{D'}(x_3')$. And so, Claim \ref{eq-lwz4} leads to
$$
|x_1'-x_2'|< 2e^{4b^2M+C}\delta_{D'}(x_3'),
$$ as required.\medskip

Based on Claim \ref{eq-lwz1}, we have
\bcl\label{sun-1}
$|x'_1-x'_3|\leq 2\theta'(\mu)e^{4b^2M+C}\delta_{D'}(x'_3).$
\ecl

In order to exploit Lemma \Ref{lem-ll-0} to show this claim, we need some preparation. It follows from \eqref{sat-5} and \eqref{sun-2} that
\begin{eqnarray*} k_{D'}(x'_1,x'_3) &\geq& \ell_{k_{D'}}(\gamma'[x'_1,x'_3])-\varepsilon
\geq \log\Big(1+\frac{\ell(\gamma'[x'_1,x'_3])}{\delta_{D'}(x'_1)}\Big)-1
\\ \nonumber&=&  \log\big(1+\frac{A}{2}\big)-1.
\end{eqnarray*} Hence, by Lemma \Ref{BHK-lem}, we have
\begin{eqnarray*}\log\left(1+\frac{|x_1-x_3|}{\min\{\delta_D(x_1),\delta_D(x_3)\}}\right)&\geq&\frac{1}{4b^2}k_D(x_1,x_3)\geq \frac{1}{4b^2M}(k_{D'}(x'_1,x'_3)-C)\\&\geq&
\frac{1}{4b^2M}\Big(\log\big(1+\frac{A}{2}\big)-1-C\Big)\\ \nonumber
&>&\log(1+A_1),\end{eqnarray*}
and so \be\label{z-006}|x_1-x_3|>A_1\min\{\delta_D(x_1),\delta_D(x_3)\}>\frac{A_1}{1+\mu}\delta_D(x_3),\ee since \eqref{new-eq-4} implies
$$\delta_D(x_3)\leq \delta_D(x_1)+|x_1-x_3|\leq (1+\mu)\delta_D(x_1).$$
Again, by \eqref{new-eq-4}, we know
$$|x_1-x_3|\leq\mu\delta_D(x_1)= 2\mu|x_1-x_2|.$$

Now, we are ready to apply Lemma \Ref{lem-ll-0} to the points $x_1$, $x_2$ and $x_3$ in $D$. Since $f$ is weakly $H$-quasisymmetric and $D$ is $b$-uniform, by considering the restriction $f|_D$ of $f$ onto $D'$, we know from Lemma \Ref{lem-ll-0} that there is an increasing function $\theta'=\theta'_{b,H}$ such that
$$|x'_1-x'_3|\leq \theta'(2\mu)|x'_1-x'_2|,$$ and thus, Claim
\ref{eq-lwz1} assures that
$$
|x'_1-x'_3|\leq  2\theta'(2\mu)e^{4b^2M+C}\delta_{D'}(x'_3),
$$
which completes the proof of Claim \ref{sun-1}.
\medskip

Let us proceed the proof. To get a contradiction to the contrary assumption \eqref{main-eq-4}, we choose $x'_4\in  D'$ such that
\be\label{sun-3} |x'_3-x'_4|=\frac{1}{3c}\delta_{D'}(x'_3).\ee
Then Lemma \Ref{ll-11} implies that \begin{eqnarray*}\log\left(1+\frac{|x_3-x_4|}{\delta_{D}(x_3)}\right)&\leq& k_{D}(x_3,x_4)\leq Mk_{D'}(x_3',x_4')+C\\ \nonumber&\leq& 3cM\frac{|x_3'-x_4'|}{\delta_{D'}(x_3')}+C\leq M+C,\end{eqnarray*} which yields that \be\label{sun-3-1}|x_3-x_4|< e^{M+C}\delta_{D}(x_3).\ee
Meanwhile, Claim \ref{sun-1} and \eqref{sun-3} imply that  $$|x_1'-x_3'|\leq 2\theta'(2\mu)e^{4b^2M+C}\delta_{D'}(x_3')= 6c\theta'(2\mu)e^{4b^2M+C}|x_3'-x_4'|.$$
Now, we apply Lemma \Ref{lem-ll-0} to the points $x_1'$, $x_3'$ and $x_4'$ in $G'$.
Since $f^{-1}:$ $G'\to G$ is weakly $H$-quasisymmetric and $G'$ is $c_1$-quasiconvex,  we know from Lemma \Ref{lem-ll-0} that there is an increasing function $\theta''=\theta''_{c_1,H}$ such that
$$|x_1-x_3|\leq \theta''\Big(6c\theta'(2\mu)e^{4b^2M+C}\Big)|x_3-x_4|,$$
which, together with \eqref{z-006} and \eqref{sun-3-1}, shows that
\begin{eqnarray*}
|x_1-x_3|&\leq&  e^{M+C}\theta''\Big(6c\theta'(2\mu)e^{4b^2M+C}\Big)\delta_D(x_3)\\ \nonumber&\leq& \frac{1+\mu}{A_1}e^{M+C}\theta''\Big(6c\theta'(2\mu)e^{4b^2M+C}\Big)|x_1-x_3|\\ \nonumber&=&
\frac{1}{2}|x_1-x_3|.
\end{eqnarray*}
This obvious contradiction shows that \eqref{main-eq-1} is true.
\qed

%%%%%%%%%%%%%%%%%
\subsection{The proof of \eqref{main-eq-2}}
 Let
$$B=12cA^2e^{6b^2M\mu\big(1+\theta''\big(\frac{1+12cA}{3c}\big)\big)},$$
and suppose on the contrary that \be\label{eq-lwz3}\ell(\gamma')> B|x'-y'|.\ee
Since $\frac{9}{2}ce^{\frac{3}{2}}<B$, we see from Lemma \ref{mon-4} that
\be \label{eq-lwz2}|x'-y'|>\frac{1}{3c} \max\{\delta_{D'}(x'), \delta_{D'}(y')\}.\ee
For convenience, in the following, we assume that $$\max\{\delta_{D'}(x'), \delta_{D'}(y')\}=\delta_{D'}(x').$$

First, we choose some special points from $\gamma'$.
By \eqref{eq-lwz3}, we know that there exist $w'_1$ and $w'_2\in \gamma'$ such that $x'$, $w'_1$, $w'_2$ and $y'$ are successive points in $\gamma'$ and
\be\label{eq-11-1}\ell(\gamma'[x',w'_1])=\ell(\gamma'[w'_2,y'])=6cA|x'-y'|.\ee
Then we have
\bcl\label{zzz-002} $|x'-w'_1|\geq \frac{1}{2}\delta_{D'}(w'_1)$ and $|y'-w'_2|\geq \frac{1}{2}\delta_{D'}(w'_2)$.
\ecl

Obviously, it suffices to show the first inequality in the claim. Suppose
 $$|x'-w'_1|< \frac{1}{2}\delta_{D'}(w'_1).$$ Then \eqref{main-eq-1} and \eqref{eq-lwz2} lead to
$$\delta_{D'}(x')\geq \delta_{D'}(w'_1)-|x'-w'_1|>\frac{1}{2}\delta_{D'}(w'_1)\geq \frac{1}{2A}\ell(\gamma'[x',w'_1])=3c|x'-y'|>\delta_{D'}(x').$$
This obvious contradiction completes the proof of Claim \ref{zzz-002}.\medskip

By using Claim \ref{zzz-002}, we get a lower bound for $|w_1-w_2|$ in terms of $\min\{\delta_D(w_1),\; \delta_D(w_2)\}$, which is as follows.
\bcl\label{eq-l2}
$|w_1-w_2|> \Big(1+\theta''\Big(\frac{1+12cA}{3c}\Big)\Big)\mu\min\{\delta_D(w_1),\; \delta_D(w_2)\}.$
\ecl
Without loss of generality, we assume that $\min\{\delta_D(w_1),\; \delta_D(w_2)\}=\delta_{D}(w_1)$.
Then by \eqref{eq-11-1} and Claim \ref{zzz-002}, we have
\be\label{eq-l10}\delta_{D'}(w'_1)\leq 2|x'-w'_1|\leq 2\ell(\gamma'[x',w'_1])=12cA|x'-y'|.\ee
Since $\gamma'$ is an $\varepsilon$-short arc and $D$ is $b$-uniform, by Lemma \Ref{BHK-lem}, we have
\begin{eqnarray*}
\log\left(1+\frac{|w_1-w_2|}{\delta_D(w_1)}\right) &\geq& \frac{1}{4b^2}k_D(w_1,w_2) \geq \frac{1}{4Mb^2}k_{D'}(w'_1,w'_2)-\frac{C}{4Mb^2}\\&\geq& \frac{1}{4Mb^2}\ell_{k_{D'}}(\gamma'[w'_1,w'_2])-\frac{\varepsilon+C}{4Mb^2}
\\ \nonumber&\geq& \frac{1}{4Mb^2}\log\left(1+\frac{\ell(\gamma'[w'_1,w'_2])}{\delta_{D'}(w'_1)}\right)-\frac{1+C}{4Mb^2}
\\ \nonumber&\geq&  \frac{1}{4Mb^2}\log\Big(1+\frac{B-12cA}{12cA}\Big)-\frac{1+C}{4Mb^2}
\\ \nonumber &=&\lambda,
\end{eqnarray*}
where the last inequality follows from \eqref{eq-l10} and the following inequalities:
$$\ell(\gamma'[w_1',w_2'])=\ell(\gamma')-\ell(\gamma'[x',w_1'])-\ell(\gamma'[y',w_2'])>(B-12cA)|x'-y'|.$$ Hence $$|w_1-w_2|\geq(e^{\lambda}-1)\delta_D(w_1)> \Big(1+\theta''\Big(\frac{1+12cA}{3c}\Big)\Big)\mu\delta_D(w_1),$$ as required. \medskip

Next, we get the following upper bound for $|w_1-w_2|$ in terms of $\min\{\delta_D(w_1),\; \delta_D(w_2)\}$.
\bcl\label{mon-3}
$|w_1-w_2|\leq \theta''\left(\frac{1+12cA}{3c}\right)\mu \min\{\delta_D(w_1),\; \delta_D(w_2)\}.$
\ecl
First, we see that $w_0\in\gamma[w_1,y]$, where $w_0$ is the point in $\gamma$ which satisfies \eqref{wes-1}, because otherwise \eqref{new-eq-4} gives that $$|w_1-w_2|\leq\mu\delta_D(w_1),$$ which contradicts with Claim \ref{eq-l2}.

%\begin{figure}[htbp]
%\begin{center}
%\input{figure04.pspdftex}
%\caption{The points $w_0$, $w_1$ and $w_2$ in $\gamma$.}
%
%\label{subinv-fig04}
%\end{center}
%\end{figure}

%{\color{red}Figure $4$: The points $w_1$, $w_2$ and $w_0$ in $\gamma$.}

We are going to apply Lemma \Ref{lem-ll-0} to the points $x'$, $w_1'$ and $w_2'$ in $G'$. We need a relationship between $|w'_1-w'_2|$ and $|x'-w_1'|$.
To this end, it follows from \eqref{eq-11-1} that $$|w'_1-w'_2|\leq |w'_1-x'|+|x'-y'|+|y'-w'_2|\leq (1+12cA)|x'-y'|\leq\frac{1+12cA}{3c}|x'-w_1'|,$$
since we infer from
the choice of $w'_1$, \eqref{main-eq-1} and Claim \ref{zzz-002} that
$$|x'-w_1'|\geq \frac{1}{2}\delta_{D'}(w_1')\geq \frac{1}{2A}\ell(\gamma'[x',w'_1])=3c|x'-y'|.$$
Then by Lemma  \Ref{lem-ll-0}, we have known that there is an increasing function $\theta''=\theta''_{c_1,H}$ such that
 $$|w_1-w_2|\leq \theta''\left(\frac{1+12cA}{3c}\right)|x-w_1|,$$
and thus, \eqref{new-eq-4} leads to
$$|w_1-w_2|\leq  \theta''\left(\frac{1+12cA}{3c}\right)\mu \min\{\delta_D(w_1),\; \delta_D(w_2)\},$$
which shows that Claim \ref{mon-3} holds.\medskip

It follows from Claims \ref{eq-l2} and \ref{mon-3} that it is impossible, and so this obvious contradiction completes the proof of \eqref{main-eq-2}.\qed

Inequalities \eqref{main-eq-1} and \eqref{main-eq-2}, together with the arbitrariness of the choice of $x'$ and $y'$ in $D'$, show that $D'$ is $B$-uniform, which implies that Theorem \ref{thm1.1} holds.\qed

\section{The proof of Theorem \ref{thm1.2} }\label{sec-5}
In this section, we always assume that $X$ and $Y$ are both $c$-quasiconvex and complete  metric spaces, that
$G \varsubsetneq X$ is a non-point-cut and locally $a$-John domain, and that $G'\varsubsetneq Y$ is a $b_1$-uniform and locally $a$-John domain. Further, we assume that $f:$ $G\to G'$ is a $\varphi$-FQC mapping and $D$ is a $b_2$-uniform subdomain of $G$.

We divide this section into two subsections. In the first  subsection, a useful lemma will be proved, and the proof of Theorem \ref{thm1.2} will be presented in the second subsection.

\subsection{An auxiliary result}

First, based on \cite[Theorems 1 and 2]{HWZ}, we prove the following result which plays a key role in the proof of Theorem \ref{thm1.2}.

\begin{lem}\label{LHWZ-lem2} Under the given assumptions in the first paragraph of this section, we have the following assertions.
\begin{enumerate}
\item\label{mon-1}
 There exist constants $M=M(c,\varphi)\geq 1$ and $C=C(c,\varphi)\geq0$ such that both $f:$ $G\to G'$
and $f^{-1}:$ $G'\to G$ are fully $(M,C)$-CQH.

\item\label{mon-2}
 There exists a constant $q=q(c)\in (0,1)$ such that for any $x\in G$, $f$ is $\eta$-quasisymmetric in $\mathbb{B}(x, q\delta_G(x))$ and $f^{-1}$ is $\eta$-quasisymmetric in $\mathbb{B}(x', q\delta_{G'}(x'))$, where $\eta=\eta_{a,c,\varphi}$.
 \end{enumerate}
\end{lem}

\bpf
Since $f$ is a $\varphi$-FQC mapping, we see that Lemma \ref{LHWZ-lem2}\eqref{mon-1} easily follows from \cite[Theorem 1]{HWZ}. For the proof of the second assertion, we infer from  \cite[Theorem 2]{HWZ} that we only need to prove that there is a homeomorphism $\theta:[0,1)\to [0, \infty)$  such that both $f$ and $f^{-1}$ are $\theta$-relative.
 By symmetry, we know that we only need to show that $f$ is $\theta$-relative. To this end, we let $0<t<1$ and $x, y\in G$ with $|x-y|=t\delta_G(x)$. Then
we separate the proof into two cases. For the first case, that is,
  $0<t \leq \frac{1}{3c}$, it follows from Lemma \Ref{ll-11} that $$k_G(x,y)\leq 3c\frac{|x-y|}{\delta_G(x)}=3c t.$$ Hence, \beq\label{lwz-eq-1}\frac{|x'-y'|}{\delta_{G'}(x')}\leq e^{k_{G'}(x',y')}-1\leq e^{\varphi(k_{G}(x,y))}-1\leq e^{\varphi(3ct)}-1.\eeq

For the other case, that is, $\frac{1}{3c}<t <1$, by Lemma \ref{LHWZ-lem1}, we know that $$k_G(x,y)\leq 2a\frac{3+t}{1-t},$$ which implies \beq\label{lwz-eq-2}\frac{|x'-y'|}{\delta_{G'}(x')}\leq e^{k_{G'}(x',y')}-1\leq e^{\varphi(k_{G}(x,y))}-1\leq e^{\varphi\big(2a\frac{3+t}{1-t}\big)}-1.\eeq
Therefore, \eqref{lwz-eq-1} and \eqref{lwz-eq-2} show that $f$ is $\theta$-relative, where

$$\theta(t)=\begin{cases}
\displaystyle e^{\varphi\big(\frac{6ac(9c+1)}{3c-1}t\big)}-1,\;
\mbox{if}\;\;t\in[0,\frac{1}{3c}],\\
\displaystyle \;\;\;e^{\varphi\big(2a\frac{3+t}{1-t})}-1,\;\;\;\, \mbox{if}\;\; t\in(\frac{1}{3c}, 1).
\end{cases}
$$
\medskip

%In the following, we show that $f^{-1}$ is $\theta_2$-relative. Let $0<t<1$ and $x', y'\in G'$ with $|x'-y'|=t\delta_{G'}(x')$. Then $$\delta_{G'}(y')\geq \delta_{G'}(x')-|x'-y'|=(1-t)\delta_{G'}(x').$$
%Since $G'$ is a $b_1$-uniform domain, we infer from Lemma \Ref{BHK-lem} that $$k_{G'}(x',y')\leq 4b_1^2\log\left(1+\frac{|x'-y'|}{\min\{\delta_{G'}(x'),\delta_{G'}(y')\}}\right)\leq 4b_1^2\log\frac{1}{1-t}, $$ which shows that $$\frac{|x-y|}{\delta_G(x)}\leq e^{k_{G}(x,y)}-1\leq  e^{\varphi(k_{G'}(x',y'))}-1\leq e^{\varphi\left(4b_1^2\log\frac{1}{1-t}\right)}-1,$$  and thus $f^{-1}$ is $\theta_2$-relative, where $$\theta_2(t)=e^{\varphi\left(4b_1^2\log\frac{1}{1-t}\right)}-1$$ for $t\in[0,1).$
Hence the proof of Lemma \ref{LHWZ-lem2} is complete.\epf

%%%%%%%%%%%%%%%%%%%%%%%%%%%%%%%%%%%%%%%%%%%%%%%%%%
%%%%%%%%%%%%%%%%%%%%%%%%%%%%%%%%%%%%%%%%%%%%%%%%%%%%%%%%%%%%%
\subsection{The proof of Theorem \ref{thm1.2}}
 For $x', y' \in D'$,  there is an $\varepsilon$-short arc $\gamma'$ joining $x'$ and $y'$ in $D'$ with
 $$0<\varepsilon<\frac{1}{2}\min\{2,k_{D'}(x',y')\}.$$
Lemmas \ref{ll-001} and \ref{LHWZ-lem2}\eqref{mon-1} show that the preimage $\gamma$ of $\gamma'$ is a $(\nu,h)$-solid arc in $D$,
 where $\nu=\nu(c, \varphi)$, $h=h(c, \varphi)$. Let $z_0\in\gamma$ be such that
$$\delta_D(z_0)=\max_{p\in \gamma}\delta_D(p).$$

%\begin{figure}[htbp]
%\begin{center}
%\input{figure06.pspdftex}
%\caption{The arc $\gamma$ and the related points.}
%
%\label{subinv-fig06}
%\end{center}
%\end{figure}

\noindent Then by Lemma \ref{ll-15}, there is a constant $\mu=\mu(b_2,\nu,h)=\mu(b_2,c,\varphi)>1$ such that for each $u\in\gamma[x,z_0]$ and for all $z\in\gamma[u,z_0]$, \be\label{new-eq-6}|u-z|\leq\mu\delta_D(z),\ee
and for each $v\in\gamma[y,z_0]$ and for all $z\in\gamma[v,z_0]$, \be\label{new-eq-7}|v-z|\leq\mu\delta_D(z).\ee

In the following, we show that $\gamma'$ is a double cone arc in $D'$.
 The proof is divided into two steps which are given in the following two subsubsections.

\subsubsection{The verification of the cigar condition} To this end, we let $$\lambda_1=\lambda_2^2(1+\lambda^4_{2})^{9(C+1)M^2b_2^2\lambda_2}$$ and $$\lambda_2=\max\Big\{\frac{3c\mu }{\varphi^{-1}(\frac{q}{21ce})},\frac{12c(C+1)\mu}{q},16b_1^2b_2^2(C+1)M^2,\frac{2(1+\mu)}{\eta^{-1}(\frac{1}{3}\eta^{-1}(\frac{1}{\mu}))}\Big\},$$
where the constants $C$, $M$, $q$ and the function $\eta$ are from Lemma \ref{LHWZ-lem2}.

In this subsubsection, we show that $\gamma'$ satisfies the cigar condition with constant $2\lambda_1^2$.
Obviously, we may assume that \begin{enumerate}
 \item\label{sun-4} either there exists a point $z_1'\in \gamma'[x',z'_0]$  such that $$\ell(\gamma'[x',z_1'])>\lambda_{1}\delta_{D'}(z_1'),$$
\item\label{sun-5}
or there exists a point  $z_2'\in \gamma'[y',z'_0]$ such that $$\ell(\gamma'[y',z_2'])>\lambda_{1}\delta_{D'}(z_2').$$
\end{enumerate}

 If \eqref{sun-4} happens, then we let $x'_0\in\gamma'[x',z'_0]$ be the first point from $x'$ to $z'_0$ such that \be\label{eq-lwz5}\ell(\gamma'[x',x'_0])=\lambda_{1}\delta_{D'}(x'_0).\ee

%\begin{figure}[htbp]
%\begin{center}
%\input{figure05.pspdftex}
%\caption{The arc $\gamma'$ and the related points.}
%
%\label{subinv-fig05}
%\end{center}
%\end{figure}

%{\color{red}Figure $5$: The arc $\gamma'$ and the related points.}

\noindent The following comparison result is useful.

\bcl\label{z-01} For every $z\in\gamma[x_0,z_0]$, $\delta_G(z)\leq\lambda_2\delta_D(z)$. \ecl

Suppose on the contrary that there exists some point $x_1\in\gamma[x_0,z_0]$ with
\be\label{mon-5}\delta_G(x_1)>\lambda_2\delta_D(x_1).\ee

%{\color{red} Figure 6: The arc $\gamma$ and the related points.}

 To get a contradiction, we need some preparation.
  First, it follows from \eqref{new-eq-6} that for all $z\in\gamma[x,x_1]$,
\be\nonumber |z-x_1|\leq\mu\delta_D(x_1)<\frac{\mu}{\lambda_2}\delta_G(x_1)\leq\frac{q}{12c(C+1)}\delta_G(x_1),\ee
whence
\be\label{fri-1}
\gamma[x, x_1]\subset \mathbb{B}\big(x_1, \frac{\mu}{\lambda_2}\delta_G(x_1)\big)\subset \mathbb{B}\big(x_1, \frac{q}{12c(C+1)}\delta_G(x_1)\big).
\ee

Further, we show that
\be\label{fri-2}
\gamma'[x', x'_1]\subset \mathbb{B}\big(x'_0, \frac{q}{10c}\delta_{G'}(x'_0)\big).
\ee

For each $z'\in\gamma'[x', x'_1]$, by Lemma \Ref{ll-11}, we have
$$k_{G'}(z',x'_1)\leq\varphi(k_G(z,x_1))\leq \varphi\big(3c\frac{|x_1-z|}{\delta_G(x_1)}\big)\leq\varphi\big(\frac{3c\mu}{\lambda_2}\big)\leq\frac{q}{21ce},$$
and so by the elementary inequality ``$e^x-1\leq ex$" in $(0,1)$, we get that $$|z'-x'_1|\leq(e^{k_{G'}(z',x'_1)}-1)\delta_{G'}(x'_1)\leq \frac{q}{21c}\delta_{G'}(x'_1),$$
which implies $$\delta_{G'}(x'_0)\geq\delta_{G'}(x'_1)-|x'_0-x'_1|\geq\frac{19}{20c}\delta_{G'}(x'_1).$$
Hence
$$ |z'-x'_0|\leq|z'-x'_1|+|x'_0-x'_1|\leq \frac{2q}{21c}\delta_{G'}(x'_1)< \frac{q}{10c}\delta_{G'}(x'_0),$$ as required.

Next, we need to choose some special points. Let $x'_3\in\gamma'[x',x'_0]$ be such that $$\ell(\gamma'[x',x'_3])=\frac{1}{2}\ell(\gamma'[x',x'_0]),$$
and let $x'_2\in\partial D'$ satisfy
$$|x'_2-x'_0|\leq \frac{3}{2}\delta_{D'}(x'_0).$$
Then we assert that
\be\label{fri-3}
x'_2,\; x'_3\in \mathbb{B}(x'_0,q\delta_{G'}(x'_0)).
\ee

It follows from \eqref{fri-2} that we only need to verify the truth of $x'_2\in \mathbb{B}(x'_0,q\delta_{G'}(x'_0))$.
Obviously, \eqref{eq-lwz5} shows that
\be\label{eq-lwz7} |x_2'-x_0'|\leq \frac{3}{2}\delta_{D'}(x_0')= \frac{3}{\lambda_1}\ell(\gamma'[x',x'_3])
< 3\delta_{D'}(x_3'),\ee
since the choice of $x'_0$ implies that $\ell(\gamma'[x',x'_3])< \lambda_1\delta_{D'}(x_3')$.
Also, we have
\begin{eqnarray}\label{eq-ll1} k_{D'}(x'_0,x'_3) &\geq& \ell_{k_{D'}}(\gamma'[x'_0,x'_3])-\varepsilon
\geq \log\left(1+\frac{\ell(\gamma'[x'_0,x'_3])}{\delta_{D'}(x'_0)}\right)-1
\\ \nonumber&=&  \log\big(1+\frac{\lambda_1}{2}\big)-1>1,
\end{eqnarray}
whence  \be\label{eq-lwz6}\nonumber |x_0'-x_3'|> \frac{1}{3c}\delta_{D'}(x_0'),\ee because otherwise, by Lemma \Ref{ll-11}, we get $$k_{D'}(x_0',x_3')< \frac{1}{3c}\frac{|x_0'-x_3'|}{\delta_{D'}(x_0')}\leq \frac{1}{9c^2}< 1,$$ which contradicts with \eqref{eq-ll1}.
Hence we deduce from \eqref{fri-2} and the choice of $x_2'$ that
\be\label{z-06} |x'_2-x'_0|\leq\frac{3}{2}\delta_{D'}(x'_0)\leq \frac{9c}{2} |x_0'-x_3'|\leq\frac{9 q}{20}\delta_{G'}(x'_0),\ee
as needed.\medskip

It is \eqref{fri-3} that allows us to apply Lemma \ref{LHWZ-lem2}\eqref{mon-2} to the points $x_0'$, $x_2'$ and $x_3'$, which shows
$$\frac{|x_2-x_0|}{|x_3-x_0|}\leq \eta\left(\frac{|x'_2-x'_0|}{|x'_3-x'_0|}\right),$$
 and thus, \eqref{eq-lwz7} implies
\be\label{z-071} |x'_3-x'_0|\leq \frac{1}{\eta^{-1}(\frac{1}{\mu})}|x'_2-x'_0|< \frac{3}{\eta^{-1}(\frac{1}{\mu})}\delta_{D'}(x'_3),\ee
since \eqref{new-eq-6} and the choice of $x_2'$ assure that $$|x_3-x_0|\leq\mu\delta_D(x_0)\leq\mu|x_2-x_0|.$$

Still, we need an estimate on $|x_3-x_0|$. Since it follows from Lemma \Ref{BHK-lem} and \eqref{eq-ll1} that
\begin{eqnarray*}\log\Big(1+\frac{|x_3-x_0|}{\min\{\delta_D(x_0),\delta_D(x_3)\}}\Big)&\geq&\frac{1}{4b_2^2}k_D(x_0,x_3)\geq \frac{1}{4b_2^2M}(k_{D'}(x'_0,x'_3)-C)\\ &\geq&
\frac{1}{4b_2^2M}\Big(\log\big(1+\frac{\lambda_1}{2}\big)-1-C\Big)\\ &>&\log(1+\lambda_2),\end{eqnarray*}
we know  \be\label{eq-ll2}|x_3-x_0|>\lambda_2\min\{\delta_D(x_0),\delta_D(x_3)\}>\frac{\lambda_2}{1+\mu}\delta_D(x_3),\ee since by \eqref{new-eq-6}, $\delta_D(x_3)\leq \delta_D(x_0)+|x_3-x_0|\leq (1+\mu)\delta_D(x_0). $\medskip

Now, we are ready to get a contradiction to the contrary assumption \eqref{mon-5}. Let $x_4\in\partial D$ be such that
$$|x_4-x_3|\leq 2\delta_D(x_3).$$
Then again by \eqref{new-eq-6} and the contrary assumption \eqref{mon-5}, we know that
\begin{eqnarray*}
|x_4-x_1| &\leq&|x_4-x_3|+|x_3-x_1|
\leq 2\delta_D(x_3)+|x_3-x_1|\leq  2\delta_D(x_1)+3|x_3-x_1|
\\ \nonumber&\leq& (2+3\mu)\delta_D(x_1)<\frac{2+3\mu}{\lambda_2}\delta_G(x_1)
\\ \nonumber&<& \frac{q}{2}\delta_G(x_1),
\end{eqnarray*}
which, together with \eqref{fri-1}, implies $$x_0,\; x_3,\; x_4\in \mathbb{B}(x_1,\frac{q}{2}\delta_G(x_1)).$$
Apply Lemma \ref{LHWZ-lem2}\eqref{mon-2} to the points $x_0,$ $x_3,$ and $x_4$. Then we see from \eqref{z-071}, \eqref{eq-ll2} and the choice of $x_4$ that
$$\frac{1}{3}\eta^{-1}\big(\frac{1}{\mu}\big)\leq \frac{|x'_4-x'_3|}{|x'_3-x'_0|}\leq \eta\left(\frac{|x_4-x_3|}{|x_3-x_0|}\right)< \eta\left(\frac{2(1+\mu)}{\lambda_2}\right),$$
which is the desired contradiction since $\lambda_2\geq \frac{2(1+\mu)}{\eta^{-1}\big(\frac{1}{3}\eta^{-1}\big(\frac{1}{\mu}\big)\big)}$.
Hence Claim \ref{z-01} holds.

\medskip

%%%%%%%%%%%%%%%%%

If $(2)$ happens,  then we let $y'_0$ be the first point in $\gamma'[ y',z'_0]$ from $y'$ to $z'_0$ such that  $$\gamma'[ y',y'_0]=\lambda_1\delta_{D'}(y'_0),$$ and the similar reasoning as in the proof of Claim \ref{z-01} implies
\bcl\label{zz-01} For every $z\in\gamma[ y_0,z_0]$, we have $\delta_G(z)\leq\lambda_2\delta_D(z)$. \ecl

%%%%%%%%%%%%%%%%%%%%
%%%%%%%%%%%%%%%%%%%%

In order to prove that $\gamma'$ satisfies the cigar condition with constant $2\lambda_1^2$, we only need to consider the case where
both \eqref{sun-4} and \eqref{sun-5} happen because the proofs for other cases are similar, and in fact, the corresponding discussions are simpler. First, we partition the part $\gamma'[x_0', y_0']$ of $\gamma'$ as follows.
Let $u'_0\in\gamma'[x'_0,y'_0]$ be such that $$\delta_{D'}(u'_0)=\max_{p'\in\gamma'[x'_0,y'_0]}\delta_{D'}(p'). $$ Obviously, there exists a unique integer $k\geq0$ with $$2^k\delta_{D'}(x'_0)\leq\delta_{D'}(u'_0)<2^{k+1}\delta_{D'}(x'_0).$$

Let $v'_0=x'_0$. If $k=0$, then we let $v'_1=u'_0$. If $k>1$, then for each $i\in\{1,$ $\ldots,$ $k\}$, we let $v'_i$ be the first point in $\gamma'[x'_0,u'_0]$ from $v'_{i-1}$ to $u'_0$ such that   $$\delta_{D'}(v'_i)=2^i\delta_{D'}(x'_0),$$ and let $v'_{k+1}=u'_0$. It is possible that $v'_{k+1}=v'_{k}$. This possibility happens if and only if $\delta_{D'}(u'_0)=2^k\delta_{D'}(x'_0)$.

%\begin{figure}[htbp]
%\begin{center}
%\input{figure07.pspdftex}
%\caption{A partition of $\gamma'[x_0',u_0']$.}
%
%\label{subinv-fig07}
%\end{center}
%\end{figure}
%{\color{red}Figure 7: The partition of $\gamma'[x_0', u_0']$.}

As for this partition, we have the following assertion.

\bcl\label{zz-02} For each $i\in\{0,$ $\ldots,$ $k\}$, we have \begin{enumerate}
\item\label{sun-6}
$\ell(\gamma'[v'_i,v'_{i+1}])\leq \lambda^4_2\delta_{D'}(v'_i),$ and
\item\label{sun-7}
for  $z'\in\gamma'[v'_i,v'_{i+1}]$,
$$\delta_{D'}(v'_i)< (1+\lambda_2^4)^{20b_1^2b_2^2M^2\lambda_2}\delta_{D'}(z').$$
\end{enumerate}
\ecl
We first prove \eqref{sun-6}.
It follows from Lemmas \Ref{BHK-lem} and \ref{LHWZ-lem2}\eqref{mon-1} that
\begin{eqnarray*}\label{eq-lwz8}
\ell_{k_{D'}}(\gamma'[v'_i,v'_{i+1}]) &\leq& k_{D'}(v'_i,v'_{i+1})+\varepsilon \leq Mk_D(v_i,v_{i+1})+C+1
\\ \nonumber &\leq& 4b_2^2M\log\left(1+\frac{|v_i-v_{i+1}|}{\min\{\delta_D(v_i),\delta_D(v_{i+1})\}}\right)+C+1,
\end{eqnarray*}
and thus, Claims \ref{z-01} and \ref{zz-01}, together with Lemma \ref{LHWZ-lem2}\eqref{mon-1}, lead to
\begin{eqnarray*}
\ell_{k_{D'}}(\gamma'[v'_i,v'_{i+1}]) &\leq&
 4b_2^2 M\log\left(1+\frac{\lambda_2|v_i-v_{i+1}|}{\min\{\delta_G(v_i),\delta_G(v_{i+1})\}}\right)+C+1
\\ \nonumber&\leq& 4b_2^2 M\lambda_2k_G(v_i,v_{i+1})+C+1
\\ \nonumber&\leq& 4b_2^2 M^2\lambda_2k_{G'}(v'_i,v'_{i+1})+4b_2^2 CM\lambda_2+C+1,
\end{eqnarray*}
and finally, with the aid of Lemma \Ref{BHK-lem}, we obtain
\begin{eqnarray}\label{eq-lwz8}
\ell_{k_{D'}}(\gamma'[v'_i,v'_{i+1}]) &\leq& 16b_1^2b_2^2M^2\lambda_2\log\left(1+\frac{|v'_i-v'_{i+1}|}{\min\{\delta_{G'}(v'_i),\delta_{G'}(v'_{i+1})\}}\right)\\\nonumber&&+4b_2^2 CM\lambda_2+C+1
\\ \nonumber&<& 16b_1^2b_2^2M^2\lambda_2\log\left(1+\frac{\ell(\gamma'[v'_i,v'_{i+1}])}{\delta_{D'}(v'_i)}\right)+4b_2^2 CM\lambda_2\\ \nonumber&&+C+1.
\end{eqnarray}
Meanwhile, the choice of $v'_i$ gives $$\ell_{k_{D'}}(\gamma'[v'_i,v'_{i+1}])=\int_{\gamma'[v'_i,v'_{i+1}]}\frac{|dz|}{\delta_{D'}(z)}\geq \frac{\ell(\gamma'[v'_i,v'_{i+1}])}{2\delta_{D'}(v'_i)},$$ and so,
$$\ell(\gamma'[v'_i,v'_{i+1}])\leq \lambda^4_2\delta_{D'}(v'_i).$$ Thus Claim \ref{zz-02}\eqref{sun-6} holds.

Moreover, by \eqref{eq-lwz8} and Claim \ref{zz-02}\eqref{sun-6}, we see that for all $z'\in \gamma'[v'_i,v'_{i+1}]$,
\begin{eqnarray*}
\log\frac{\delta_{D'}(v'_i)}{\delta_{D'}(z')} &\leq& k_{D'}(v'_i,z') \leq \ell_{k_{D'}}(\gamma'[v'_i,v'_{i+1}])
\\ \nonumber&\leq& 16b_1^2b_2^2M^2\lambda_2\log(1+\lambda_2^4)+4b_2^2 CM\lambda_2+C+1
\\ \nonumber&\leq& 20b_1^2b_2^2M^2\lambda_2\log(1+\lambda_2^4).
\end{eqnarray*}
Hence
$$\delta_{D'}(v'_i)< (1+\lambda_2^4)^{20b_1^2b_2^2M^2\lambda_2}\delta_{D'}(z'),$$
which shows that Claim \ref{zz-02}\eqref{sun-7} is true too, and thus the claim is proved.\medskip

Now, we are ready to verify that $\gamma'$ satisfies the $2\lambda_{1}^2$-cigar condition, that is, for every $z'\in\gamma'[x',y']$,
\be\label{zz-05} \min\{\ell(\gamma'[x',z']),\ell(\gamma'[z',y'])\}\leq 2\lambda_{1}^2\delta_{D'}(z').\ee
Obviously, we only need to consider the case $z'\in\gamma'[x',u'_0]$ since the proof for the case  $z'\in\gamma'[y',u'_0]$ is similar.

If $z'\in\gamma'[x',x'_0]$, then \eqref{zz-05} immediately follows from the choice of $x'_0$.

 If $z'\in\gamma'[x'_0,u'_0]$, then there must exist a $j\in\{0,\ldots,k\}$ such that $$z'\in\gamma'[v'_j,v'_{j+1}],$$ and so by \eqref{eq-lwz5} and Claim \ref{zz-02}, we have the following:
\begin{eqnarray*}
\ell(\gamma'[x',z']) &=&  \ell(\gamma'[x',x'_0])+\ell(\gamma'[x'_0,z'])
\leq \lambda_1\delta_{D'}(x'_0)+\sum_{i=0}^j\ell(\gamma'[v'_i,v'_{i+1}])
\\ \nonumber&\leq& \lambda_1\delta_{D'}(v'_j)+\lambda_2^4\sum_{i=0}^j\delta_{D'}(v'_i)
\leq (\lambda_1+2\lambda_2^4)\delta_{D'}(v'_j) \\ \nonumber&<& 2\lambda_{1}^2\delta_{D'}(z'),
\end{eqnarray*}
which shows that \eqref{zz-05} is also true in this case. Hence $\gamma'$ satisfies the $2\lambda_{1}^2$-cigar condition. \qed

\subsubsection{The verification of the turning condition}
We only need to prove that
\be\label{main-lem2}\ell(\gamma')\leq \lambda_0|x'-y'|,\ee where $$\lambda_0=\max\Big\{24c\lambda_{1}^2e^{1+C}e^{6Mb_2^2[1+\mu\eta(8\lambda_{1}^2+\frac{1}{3c})]},48c\lambda_{1}^2e^{\rho}\Big\}$$ and
 $$\rho=\frac{96(c+1)(b_1b_2M)^2}{q\varphi^{-1}(\log\frac{3}{2})\eta^{-1}(\frac{q}{12\lambda_{1}^2})}
 \log\Big(1+4c\lambda_{1}^2+\frac{1}{6c}\Big)+\frac{24(c+1)MCb_2^2}{q\varphi^{-1}(\log\frac{3}{2})\eta^{-1}(\frac{q}{12\lambda_{1}^2})}+C+1.$$
We prove \eqref{main-lem2} by contradiction. Suppose that
\be\label{sun-8}\ell(\gamma')> \lambda_0|x'-y'|.\ee Then Lemma \ref{mon-4} implies  $$|x'-y'|>\frac{1}{3c}\max\{\delta_{D'}(x'),\delta_{D'}(y')\},$$ and also, we know that there exist $w'_1$, $w'_2\in \gamma'$ such that
$x'$, $w'_1$, $w'_2$ and $y'$ are successive pints in $\gamma'$, and
\be\label{eq-lwz10}
\ell(\gamma'[x',w'_1])=\ell(\gamma'[w'_2,y'])=12c\lambda_{1}^2|x'-y'|.
\ee
%
%\begin{figure}[htbp]
%\begin{center}
%\input{figure08.pspdftex}
%\caption{The points $w_1'$, $w_2'$ and $w_3'$ in $\gamma'$.} \label{subinv-fig08}
%\end{center}
%\end{figure}
%{\color{red} Figure $8$: The points $w_1'$, $w_2'$ and $w'_3$ in $\gamma'$.}

As for this partition of $\gamma'$, we prove several claims.
\bcl\label{thur-1}
 $|x'-w'_1|\geq \frac{1}{2}\delta_{D'}(w'_1)\;\;{\rm and} \;\;|y'-w'_2|\geq \frac{1}{2} \delta_{D'}(w'_2).$
\ecl
The proof of this claim easily follows from a similar argument as in the proof of Claim \ref{zzz-002} with the substitution $A$ by $2\lambda_1^2$.

\bcl\label{thur-2}
$6c|x'-y'| \leq \min\{\delta_{D'}(w'_1), \; \delta_{D'}(w'_2)\} \leq \max\{\delta_{D'}(w'_1), \; \delta_{D'}(w'_2)\} \leq 24c\lambda_{1}^2|x'-y'|.$
\ecl
It is equivalent to show that $6c|x'-y'|\leq\delta_{D'}(w'_i)\leq 24c\lambda_{1}^2|x'-y'|$ for $i\in \{1,\;2\}$, which easily follows from Claim \ref{thur-1},
\eqref{zz-05} and  \eqref{eq-lwz10}.\medskip

\bcl\label{thur-3}
$|w_1-w_2|>\mu\eta\big(8\lambda_1^2+\frac{1}{3c}\big)\min\{\delta_D(w_1), \delta_D(w_2)\}.$
\ecl
For convenience, we may assume that $$\min\{\delta_D(w_1), \delta_D(w_2)\}=\delta_D(w_1).$$
Since $D$ is $b_2$-uniform, we see from Lemma \Ref{BHK-lem} that
$$
\log\left(1+\frac{|w_1-w_2|}{\delta_D(w_1)}\right) \geq \frac{1}{4b_2^2}k_D(w_1,w_2) \geq \frac{1}{4Mb_2^2}( k_{D'}(w'_1,w'_2)-C),
$$
and thus, the assumption ``$\gamma'$ being an $\varepsilon$-short arc" implies
\begin{eqnarray*}
\log\left(1+\frac{|w_1-w_2|}{\delta_D(w_1)}\right) &\geq& \frac{1}{4Mb_2^2}(\ell_{k_{D'}}(\gamma'[w'_1,w'_2])-\varepsilon-C)
\\ \nonumber&\geq& \frac{1}{4Mb_2^2}\left(\log\left(1+\frac{\ell(\gamma'[w'_1,w'_2])}{\delta_{D'}(w'_1)}\right)-1-C\right)
\\ \nonumber&\geq&  \frac{1}{4Mb_2^2}\left(\log\left(1+\frac{\lambda_0-24c\lambda_{1}^2}{24c\lambda_{1}^2}\right)-1-C\right)\\ \nonumber &=&\mu_1,
\end{eqnarray*} where the last inequality holds because of Claim \ref{thur-2} and the following chain of inequalities which are from \eqref{sun-8} and \eqref{eq-lwz10}:
\be\label{eq-lwz13}\ell(\gamma'[w_1',w_2'])=\ell(\gamma')-\ell(\gamma'[x',w_1'])-\ell(\gamma'[y',w_2'])>(\lambda_0-24c\lambda_{1}^2)|x'-y'|.\ee
Hence \be\label{eq-l22}\nonumber |w_1-w_2|\geq(e^{\mu_1}-1)\delta_D(w_1)>\mu\eta\big(8\lambda_1^2+\frac{1}{3c}\big)\delta_D(w_1),\ee as needed.
\medskip

With the aid of \eqref{new-eq-6}, the following is a direct consequence of Claim \ref{thur-3}.
\bcl\label{thur-4}
$|x-w_1|\leq \mu\delta_D(w_1).$
\ecl

For points in $\gamma'[w'_1,w'_2]$, we have the following comparison result.
\bcl\label{zzz-003} For all $w'\in\gamma'[w'_1,w'_2]$, $\delta_{G'}(w')\leq \frac{12\lambda_{1}^2}{q}\delta_{D'}(w')$. \ecl

Suppose that there exists some point $w'_3\in\gamma'[w'_1,w'_2]$ with
 \be\label{eq-116-1}\delta_{G'}(w'_3)> \frac{12\lambda_{1}^2}{q}\delta_{D'}(w'_3).\ee
Without loss of generality, we may assume that $$\min\{\ell(\gamma'[x',w'_3]),\ell(\gamma'[y',w'_3])\}=\ell(\gamma'[x',w'_3]),$$ since the proof for  the other case is similar.

By \eqref{zz-05} and \eqref{eq-116-1}, for all $w'\in\gamma'[x',w'_3]$, we have $$\ell(\gamma'[w',w'_3])\leq 2\lambda_{1}^2\delta_{D'}(w'_3)<\frac{q}{6}\delta_{G'}(w_3'),$$
and  so, \eqref{eq-lwz10} gives  \begin{eqnarray*}|w_2'-w'_3|&\leq&|w_2'-y'|+|y'-x'|+|x'-w'_3|\\ &\leq& \ell(\gamma'[y',w_2'])+\frac{1}{12c\lambda_1^2}\ell(\gamma'[x',w_1'])+\ell(\gamma'[x',w_3']) \\ &\leq& \big(2+\frac{1}{12c\lambda_1^2}\big)\ell(\gamma'[x',w'_3])<\frac{q}{2}\delta_{G'}(w_3').\end{eqnarray*}
Then we know that $$x',\; w'_1,\; w'_2\in \mathbb{B}(w'_3,\frac{q}{2}\delta_{G'}(w'_3)).$$
By applying Lemma \ref{LHWZ-lem2}\eqref{mon-2} to the points $x'$, $w'_1$ and $w'_2$, together with \eqref{eq-l22} and Claims \ref{thur-1} $\sim$ \ref{thur-4}, we know that
$$\eta\big(8\lambda_{1}^2+\frac{1}{3c}\big)< \frac{|w_2-w_1|}{|x-w_1|}\leq \eta\left(\frac{|w'_1-w'_2|}{|x'-w'_1|}\right)\leq \eta\big(8\lambda_{1}^2+\frac{1}{3c}\big),$$
since by \eqref{eq-lwz10}, $$|w'_1-w'_2|\leq |w'_2-y'|+|y'-x'|+|x'-w'_1|\leq (1+24c\lambda_{1}^2)|x'-y'|.$$
This is the desired contradiction, from which the claim follows.
\medskip

%%%%%%%%%%%%%%%%%
The next result is an analogue of Claim \ref{zzz-003} for points in $\gamma[w_1,w_2]$.

\bcl\label{zzz-004} For all $w\in\gamma[w_1,w_2]$, we have $$\delta_{G}(w)\leq \frac{6(c+1)}{q\varphi^{-1}(\log\frac{3}{2})\eta^{-1}\big(\frac{q}{12\lambda_{1}^2}\big)}\delta_{D}(w).$$\ecl

Suppose that there exists some point $u\in\gamma[w_1,w_2]$ such that
$$\delta_{G}(u)> \frac{6(c+1)}{q\varphi^{-1}(\log\frac{3}{2})\eta^{-1}\big(\frac{q}{12\lambda_{1}^2}\big)}\delta_{D}(u).$$

%\begin{figure}[htbp]
%\begin{center}
%\input{figure09.pspdftex}
%\caption{The points $w_1$, $u$ and $w_2$ in $\gamma$, $u_1$ in $G$ and $u_2$ in $\partial D$.}
%
%\label{subinv-fig09}
%\end{center}
%\end{figure}

%{\color{red}Figure $9$: The points $w_1$, $w_2$, $u$ in $\gamma$, $u_1$ in $G$  and $u_2$ in $\partial D$.}

Let
$$u_1\in \mathbb{S}\Big(u,\frac{q\varphi^{-1}(\log\frac{3}{2})}{3(c+1)}\delta_G(u)\Big).$$
Then Lemma \ref{lem-LHWZ1} guarantees that $u_1\in G$, and thus,
 Lemma \Ref{ll-11} implies
$$k_{G'}(u_1',u')\leq \varphi(k_G(u_1,u))\leq \varphi\Big(3c\frac{|u_1-u|}{\delta_G(u)}\Big)< \log\frac{3}{2},$$ which leads to
 \be\label{eq-lwz121}|u_1'-u'|\leq\frac{1}{2}\delta_{G'}(u').\ee
%We infer from Lemma \ref{lem-LHWZ1} that $u'_1\in G'$, and  \be\label{eq-lwz12}|u_1-v|=\frac{q}{2(c+1)}\delta_G(v).\ee
Let $u_2\in\partial D$ be such that  $$|u_2-u|\leq 2\delta_D(u).$$
Thus \be\label{eq-lwz11}|u_2-u|\leq 2\delta_D(u)\leq \frac{q\varphi^{-1}(\log\frac{3}{2})}{3(c+1)}\eta^{-1}\big(\frac{q}{12\lambda_{1}^2}\big) \delta_G(u)<\frac{q}{3(c+1)}\delta_G(u).\ee
Hence,  the choice of $u_1$ implies  $$u_1,\; u_2\in \mathbb{B}(u,\frac{q}{2}\delta_G(u))\cap G.$$
Apply Lemma \ref{LHWZ-lem2}\eqref{mon-2} to the points $u$, $u_1$ and $u_2$. Then
 \eqref{eq-lwz121} and \eqref{eq-lwz11} lead to $$\frac{q}{6\lambda_{1}^2}< \frac{|u'_2-u'|}{|u'_1-u'|}\leq \eta\left(\frac{|u_2-u|}{|u_1-u|}\right)<\frac{q}{12\lambda_{1}^2},$$
 since Claim \ref{zzz-003} leads to $$|u'_2-u'|\geq \delta_{D'}(u')\geq  \frac{q}{12\lambda_{1}^2}\delta_{G'}(u').$$
 This obvious contradiction completes the proof. \medskip

Now, we are ready to get a contradiction to the contrary assumption \eqref{sun-8}.
 By  Claim \ref{thur-2} and \eqref{eq-lwz13}, we have
$$\log\left(1+\frac{\ell(\gamma'[w'_1,w'_2])}{\min\{\delta_{D'}(w'_1),\delta_{D'}(w'_2)\}}\right)
\geq \log\left(1+\frac{\lambda_0-24c\lambda_{1}^2}{24c\lambda_{1}^2}\right).$$
Also, Lemma \Ref{BHK-lem} implies
\begin{eqnarray*}
 \log\left(1+\frac{\ell(\gamma'[w'_1,w'_2])}{\min\{\delta_{D'}(w'_1),\delta_{D'}(w'_2)\}}\right)
&\leq & \ell_{k_{D'}}(\gamma'[w'_1,w'_2])\leq k_{D'}(w'_1,w'_2)+\varepsilon
\\ \nonumber&\leq& Mk_D(w_1,w_2)+C+1
\\ \nonumber&\leq& 4b_2^2M\log\left(1+\frac{|w_1-w_2|}{\min\{\delta_D(w_1),\delta_D(w_2)\}}\right)+C+1.
\end{eqnarray*}
Furthermore, Claim \ref{zzz-004} gives
$$\min\{\delta_D(w_1),\delta_D(w_2)\}\geq \frac{q\varphi^{-1}(\log\frac{3}{2})\eta^{-1}\big(\frac{q}{12\lambda_{1}^2}\big)}{6(c+1)}\min\{\delta_G(w_1),\delta_G(w_2)\},$$
and Lemma \ref{LHWZ-lem2}\eqref{mon-1} leads to
$$k_G(w_1,w_2)\leq M k_{G'}(w'_1,w'_2)+C.$$
Consequently, it follows from Lemma \Ref{BHK-lem} that
\begin{eqnarray*}
\log\left(1+\frac{\lambda_0-24c\lambda_{1}^2}{24c\lambda_{1}^2}\right)
&\leq &
 \frac{24b_2^2(c+1)M^2}{q\varphi^{-1}(\log\frac{3}{2})\eta^{-1}\big(\frac{q}{12\lambda_{1}^2}\big)}k_{G'}(w'_1,w'_2)\\ \nonumber&+&\frac{24b_2^2(c+1)C M}{q\varphi^{-1}(\log\frac{3}{2})\eta^{-1}\big(\frac{q}{12\lambda_{1}^2}\big)}+C+1
\\ \nonumber&\leq& \frac{96b_1^2 b_2^2(c+1)CM^2}{q\varphi^{-1}(\log\frac{3}{2})\eta^{-1}\big(\frac{q}{12\lambda_{1}^2}\big)}\log\left(1+\frac{ |w'_1-w'_2|}{\min\{\delta_{G'}(w'_1),\delta_{G'}(w'_2)\}}\right)\\ \nonumber&&+\frac{24b_2^2(c+1)C M}{q\varphi^{-1}(\log\frac{3}{2})\eta^{-1}\big(\frac{q}{12\lambda_{1}^2}\big)}+C+1.
\end{eqnarray*}
Since \eqref{eq-lwz10} and Claim \ref{thur-2} lead to
\beq\nonumber |w'_1-w'_2|&\leq& |w'_1-x'|+|x'-y'|+|y'-w'_2|\leq (24c\lambda_1^2+1)|x'-y'|\\&\leq&  \nonumber \big(4\lambda_{1}^2+\frac{1}{6c}\big)\min\{\delta_{G'}(w'_1),\delta_{G'}(w'_2)\}.\eeq Finally, we see that
\begin{eqnarray*}
\log\left(1+\frac{\lambda_0-24c\lambda_{1}^2}{24c\lambda_{1}^2}\right)&\leq &
 \frac{96b_1^2b_2^2(c+1)CM^2}{q\eta^{-1}\big(\frac{q\varphi^{-1}(\log\frac{3}{2})}{12\lambda_{1}^2}\big)}
\log\Big(1+4\lambda_{1}^2+\frac{1}{6c}\Big)\\ &&+\frac{24b_2^2(c+1)CM}{q\varphi^{-1}(\log\frac{3}{2})\eta^{-1}\big(\frac{q}{12\lambda_{1}^2}\big)}+C+1
\\ \nonumber&=& \rho,\end{eqnarray*}
which is the desired contradiction since $\lambda_0\geq48c\lambda_{1}^2e^{\rho}$. Hence Theorem \ref{thm1.2} holds.\qed

%by \eqref{zz-05}, Lemma \Ref{BHK-lem},  Claim \ref{zzz-004}, \eqref{eq-lwz10}, \eqref{eq-l11}, \eqref{eq-l11'} and \eqref{eq-lwz13}

\section{The proof of Theorem \ref{cor-1}}\label{sec-6}
We start with two notations.
Let $f: X\to Y$ be a homeomorphism between two metric spaces, and let $x$ be a non-isolated point of $X$. We write
$$L(x,f)=\limsup_{y\to x}\frac{|y'-x'|}{|y-x|}\;\;\mbox{and}\;\;l(x,f)=\liminf_{y\to x}\frac{|y'-x'|}{|y-x|}.$$

Suppose $G$ denotes a proper subdomain in $X$. If $|\cdot|=k_G(\cdot)$, then we denote $L(x,f)$ and $l(x,f)$ by $L_{k_G}(x,f)$ and $l_{k_G}(x,f)$, respectively.

Now, we are going to show three lemmas. The first lemma is about the comparison of the quantities $L(x,f)$ and $L_{k_G}(x,f)$ (resp. $l(x,f)$ and $l_{k_G}(x,f)$).

\begin{lem}\label{lem5-1} Suppose that $X$ is $c$-quasiconvex and $Y$ is $c'$-quasiconvex, and that
$G\subsetneq X$ and $G'\subsetneq Y$ are domains. If $f: G\to G'$ is continuous, then
$$\frac{L(x,f)\delta_G(x)}{6c\delta_{G'}(x')}\leq L_{k_G}(x,f)\leq 6c'\frac{L(x,f)\delta_G(x)}{\delta_{G'}(x')}$$ and $$\frac{l(x,f)\delta_G(x)}{6c\delta_{G'}(x')}\leq l_{k_G}(x,f)\leq 6c'\frac{l(x,f)\delta_G(x)}{\delta_{G'}(x')}.$$
\end{lem}

\bpf By symmetry, it suffices to prove the first chain of inequalities in the lemma.
On the one hand, by Lemma \Ref{ll-11}, we have \begin{eqnarray*}
L_{k_G}(x,f)&=&\limsup_{y\to x}\frac{k_{G'}(x',y')}{k_G(x,y)}=\limsup_{y\to x}\left(\frac{k_{G'}(x',y')}{|x'-y'|}\frac{|x'-y'|}{|x-y|}\frac{|x-y|}{k_G(x,y)}\right)\\&\leq&L(x,f)\frac{6c'\delta_G(x)}{\delta_{G'}(x')}.
\end{eqnarray*}
 On the other hand, again by Lemma \Ref{ll-11}, we know
\begin{eqnarray*}
L(x,f)&=&\limsup_{y\to x}\frac{|x'-y'|}{|x-y|}\\&=&\limsup_{y\to x}\left(\frac{|x'-y'|}{k_{G'}(x',y')}\frac{k_{G'}(x',y')}{k_G(x,y)}\frac{k_G(x,y)}{|x-y|}\right)\\&\leq&L_{k_G}(x,f)\frac{6c\delta_{G'}(x')}{\delta_G(x)}.
\end{eqnarray*}
Hence the proof is complete.
\epf

%\begin{Lem}\label{lem5-2} Suppose that $X$ is $c$-quasiconvex and that $f:X\to Y$ is a map with $L(x,f)\leq M$ for all $x\in X$. Then $f$ is $cM$-Lipschitz.
%\end{Lem}
The next lemma is a characterization for a homeomorphism from $X$ to $Y$ to be $M$-QH in terms of $l_{k_G}(x,f)$ and $L_{k_G}(x,f)$.

\begin{lem}\label{lem5-3} Suppose $f:G\to G'$ is a homeomorphism, where $G$ denotes a proper subdomain of $X$. Then $f$ is $M$-QH if and only if $$\frac{1}{M}\leq l_{k_G}(x,f)\leq L_{k_G}(x,f)\leq M.$$
\end{lem}

\bpf  Since the spaces $(G,k_G)$ and $(G',k_{G'})$ are $\tau$-quasiconvex for all $\tau> 1$ (see \cite[Lemma 2.5]{HWZ}), we easily know from \cite[Lemma 5.5]{Vai6} that the lemma holds.
\epf

We are ready to prove the main lemma in this section.

\begin{lem}\label{lem5-4} Suppose that $X$ is $c$-quasiconvex and $Y$ is $c'$-quasiconvex, and that
$G\subsetneq X$ and $G'\subsetneq Y$ are domains. If $f: G\to G'$ is $M$-QH, then
$f$ is fully $M'$-QH with $M'=216 (c c'M)^2$. That is, $f$ is $\varphi$-FQC with $\varphi(t)=M't$.
\end{lem}
\bpf Let $D$ be an arbitrary domain in $G$. Fix a point $x$ in $D$. By symmetry and Lemma \ref{lem5-3}, it suffices to show that  \be\label{mon-10} L_{k_D}(x,f)\leq M'.\ee
Since $f$ is $M$-QH, Lemma \ref{lem5-3} gives $L_{k_G}(x,f)\leq M$. Hence, to prove \eqref{mon-10}, by Lemma \ref{lem5-1}, it suffices to show that \be\label{eq-lem5-4-1}\frac{\delta_D(x)}{\delta_{D'}(x')}\leq 6cc'M \frac{\delta_G(x)}{\delta_{G'}(x')}.\ee

If $\delta_{D'}(x')\geq \frac{1}{6cc'M}\delta_{G'}(x')$, then \eqref{eq-lem5-4-1} is obvious. So we assume in the following that $$\delta_{D'}(x')< \frac{1}{6cc'M}\delta_{G'}(x').$$ Let $$0<\epsilon<\frac{1}{6cc'M}\delta_{G'}(x')-\delta_{D'}(x'),$$ and let $y'\in \partial D'$ be such that \be\label{eq-lwz14}|y'-x'|<\delta_{D'}(x')+\epsilon<\frac{1}{6cc'M}\delta_{G'}(x').\ee  Then by Lemma \ref{lem-LHWZ1}, we know $y'\in G'$, and thus,
Lemma \Ref{ll-11} leads to $$k_{G'}(x',y')\leq 3c'\frac{|x'-y'|}{\delta_{G'}(x')}\leq\frac{1}{2cM},$$ which shows that $$k_G(x,y)\leq M k_{G'}(x',y')\leq \frac{1}{2c}.$$ Hence again by Lemma \Ref{ll-11}, together with \eqref{eq-lwz14}, we have
\begin{eqnarray*}
\delta_D(x)&\leq& |x-y|\leq 2k_G(x,y)\delta_G(x)\leq 2Mk_{G'}(x',y')\delta_G(x) \leq 6c'M\frac{|x'-y'|}{\delta_{G'}(x')}\delta_G(x)\\&\leq& 6c'M\frac{\delta_{D'}(x')+\epsilon}{\delta_{G'}(x')}\delta_G(x).
\end{eqnarray*} Letting $\epsilon\to 0$ gives \eqref{eq-lem5-4-1}. Hence the proof of Lemma \ref{lem5-4} is complete.
\epf

{\bf The proof of Theorem \ref{cor-1}}.  Obviously, the proof of Theorem \ref{cor-1} easily follows from Theorem \ref{thm1.2} and Lemma \ref{lem5-4}.
%{\color{red}???????}

\end{document}